\documentclass[a4,12pt]{article}

\usepackage{amsmath, amssymb, amsthm, multicol}

\usepackage[english]{babel}

\usepackage{xcolor}

\newtheorem{lemma}{Lemma}
\newtheorem{theorem}{Theorem}
\newtheorem{definition}{Definition}
\newtheorem{prop}{Proposition}
\newtheorem{coro}{Corollary}

\allowdisplaybreaks

\theoremstyle{remark}

\def\N{{\mathbb N}} 
\def\Z{{\mathbb Z}}
\def\Q{{\mathbb Q}}
\def\R{{\mathbb R}} 
\def\C{{\mathbb C}}

\def\s{{\bf s}}
\def\m{{\bf m}}
\def\e{{\bf e}}
\def\x{{\bf x}}
\def\y{{\bf y}}

\def\a{{\bf a}}
\def\b{{\bf b}}
\def\k{{\bf k}}

\def\X{{\bf X}}

\def\u{{\bf u}}

\def\d{{\bf d}}

\def\Nb{{\bf N}}
\def\P{{\bf P}}
\def\H{{\bf H}}

\newcommand{\zerob}{\boldsymbol{0}}
\newcommand{\unb}{\boldsymbol{1}}

\newcommand{\alphab}{\boldsymbol{\alpha}}
\newcommand{\betab}{\boldsymbol{\beta}}

\newcommand{\gammab}{\boldsymbol{\gamma}}
\newcommand{\thetab}{\boldsymbol{\theta}}

\newcommand{\nub}{\boldsymbol{\nu}}
\newcommand{\mub}{\boldsymbol{\mu}}

\def\eps{{ \varepsilon }}
\newcommand{\be}{\begin{enumerate}}
\newcommand{\ee}{\end{enumerate}}
\let\ds=\displaystyle

\begin{document}
\title{Values of multiple zeta-functions with polynomial denominators at non-positive integers\footnote{The authors benefit from the financial support of 
the French-Japanese Project ``Zeta Functions of Several Variables and  Applications" 
(PRC CNRS/JSPS 2015-2016).}}
\author{Driss Essouabri \hskip 0.4cm and \hskip 0.2cm
Kohji Matsumoto}
\date{}
\maketitle 
\vspace{10mm}

\noindent
{{\bf {Abstract.}}
We study rather general multiple zeta-functions whose denominators are given by
polynomials.     The main aim is to prove explicit formulas for the values of those
multiple zeta-functions at non-positive integer points.    We first treat the case when
the polynomials are power sums, and observe that some ``trivial zeros'' exist.
We also prove that special values are sometimes transcendental.   Then we proceed to
the general case, and show an explicit expression of special values at non-positive
integer points which involves certain period integrals.    We give examples of
transcendental values of those special values or period integrals.    We also mention
certain relations among Bernoulli numbers which can be deduced from our explicit
formulas.    Our proof of explicit formulas are based on the Euler-Maclaurin
summation formula, Mahler's theorem, and a Raabe-type lemma due to Friedman and Pereira.

\bigskip

\noindent
{\small {\bf Mathematics Subject Classifications: 11M32, 11J81} }
\setcounter{tocdepth}{2}
\vspace{10mm}

\section {Introduction}\label{sec-intro}

Denote by $\N, \N_0, \Z, \Q, \R, \C$ the set of positive integers, non-negative integers, 
rational integers, rational numbers, real numbers, and complex numbers, respectively.

Let $\gammab =(\gamma_1,\dots, \gamma_n) \in \C^n$ and  $\b =(\b_1,\dots, b_n) \in \C^n$ be two vectors of complex parameters such that 
$\Re (\gamma_j) >0$ and $\Re (b_j) >-\Re (\gamma_1)$ for all $j=1,\dots, n$.
The generalized  Euler-Zagier multiple zeta-function,
introduced in \cite{mat-JNT}, is defined  
for $n-$tuples of complex variables $\s=(s_1,\dots, s_n)$ by
\begin{equation}\label{EZMzetadef}
\zeta_n(\s; \gammab; \b) :=
\sum_{m_1\geq 1 \atop m_2, \dots, m_n \geq 0} \frac{1}{\prod_{j=1}^n (\gamma_1 m_1+\dots+ \gamma_j m_j +b_j)^{s_j}}.
\end{equation}
If $b_1=0$ and $b_j=\gamma_2+\dots+\gamma_j$ for all $j=2,\dots, n$ then $\zeta_n(\s; \gammab; \b)$ coincides with the multiple zeta-function 
$\zeta_n(\s; \gammab)$ considered in \cite{FKMT}.
If in addition $\gamma_j=1$ for all $j=1,\dots, n$, then  $\zeta_n(\s; \gammab; \b)$ coincides with the classical Euler-Zagier multiple zeta-function
 $\ds
\sum_{1\leq m_1< m_2<\dots <m_n } \frac{1}{m_1^{s_1}\dots m_n^{s_n}}.$\par
The generalized  Euler-Zagier multiple zeta-function $\zeta_n(\s; \gammab; \b)$ converges absolutely in the domain 
\begin{equation}\label{domaincv}
\mathcal D_n:=\{\s=(s_1,\dots, s_n)\in \C^n \mid \Re (s_j+\dots+ s_n) >n+1-j~~ 
(1\leq j\leq n)\},
\end{equation}
and has the meromorphic continuation to the whole complex space $\C^n$ whose poles are located in the union of the hyperplanes 
$$s_j+\dots+s_n =(n+1-j)-k_j \quad (1\leq j\leq n,~k_1,\dots, k_n \in \N_0).$$
Moreover it is known that for $n\geq 2$, almost all non-positive integer points lie on the singular locus above and are point of indeterminacy.
In \cite{komori}, Y. Komori proved that for any $\Nb =(N_1,\dots, N_n)\in \N_0^n$ and $\thetab=(\theta_1,\dots, \theta_n) \in \C^n$ such that $\theta_j+\dots+\theta_n \neq 0$ for all $j=1,\dots, n$, 
the limit 
\begin{equation}\label{gmzvtheta}
\zeta_n^{\thetab} (-\Nb; \gammab; \b):=\lim_{t\rightarrow 0}\zeta_n(-\Nb+t\thetab; \gammab; \b)
\end{equation} 
exists and  express this limit in term of $\Nb$, $\thetab$ and generalized Bernoulli numbers defined implicitly as coefficients of some multiple series. \par
In our recent work \cite{GMZV}, we proved  a closed explicit formula for $\zeta_n^{\thetab} (-\Nb; \gammab; \b)$ in terms of $\Nb$, $\thetab$ and {\it only classical Bernoulli numbers} 
$B_n$ defined by
\begin{equation}\label{def_Bernoulli}
\frac{X}{e^X-1}=\sum_{n=0}^{\infty}\frac{B_n}{n!}X^n.
\end{equation}
We also gave several results on the values of its partially twisted analogues
in \cite{Twistedcase}.

The aim of the present paper is to consider the case of multiple Dirichlet series whose
denominator is given by polynomials of any degree.
A natural non-linear extension of the Euler-Zagier multiple zeta-function is the
series defined for $n-$tuples of complex variables $\s=(s_1,\dots, s_n)$ by
\begin{equation}\label{NEZMzetadef}
\zeta_{n, \d, \gammab} (\s) :=
\sum_{m_1, \dots, m_n \geq 1} \frac{1}{\prod_{j=1}^n (\gamma_1 m_1^{d_1}+\dots+ \gamma_j m_j^{d_j} )^{s_j}}.
\end{equation}
We begin with the discussion of this type of series, because of the following two
reasons.     First, this series is ``not so far'' from the Euler-Zagier multiple
zeta-function \eqref{EZMzetadef}, so we can find several common features with 
\eqref{EZMzetadef}, or even with the Riemann zeta-function (such as {\it trivial zeros}).
Secondly, on the other hand, some propeties different from the linear case
already appear in this special type of non-linear case (such as {\it the transcendency
of special values}).
In Section \ref{sec-firstmainextension} we will state the main
results on \eqref{NEZMzetadef} (Proposition \ref{basicproperties1}, Theorem
\ref{main-main}, Theorem \ref{irrationalityTh} and corollaries).
We will prove Proposition \ref{basicproperties1} in Section \ref{sec-proofstheorems123},
Theorem \ref{main-main} and its corollaries in Section \ref{sec-proof_th1_cor123},
and then Theorem \ref{irrationalityTh} in Section \ref{sec-proofirrationalityTh}.

Then we proceed to the discussion of more general multiple series with
polynomial denominators.
Consider for any $j=1,\dots, n$ a polynomial $ P_j\in \R[X_1,\dots, X_j]$ in $j$ variables and assume that for all $j=1,\dots, n$:
\begin{equation}\label{growthcondition}
 P_j(x_1,\dots,x_j) \rightarrow \infty {\mbox { as }} x_1+\dots+x_j \rightarrow \infty, ~\left(\x=(x_1,\dots,x_j)\in [1,\infty)^j\right).
\end{equation}
We assume here for simplicity that for all $j=1,\dots, n$,
\begin{equation}\label{forsimplicity}
P_j(x_1,\dots, x_j) >0 \quad \mbox{for all}\;\; (x_1,\dots, x_j)\in [1,\infty)^j,
\end{equation}
and define the multiple zeta-function with polynomial denominators 
for $n-$tuples of complex variables $\s=(s_1,\dots, s_n)$ by
\begin{equation}\label{NLEZMzetadef}
\zeta_n(\s; \P) :=
\sum_{ m_1, \dots, m_n \geq 1} \frac{1}{\prod_{j=1}^n P_j(m_1, \dots, m_j)^{s_j}},
\end{equation}
where $\P=(P_1,\ldots,P_n)$.

By using Lemma 1 of \cite{essouabriFourier}, the condition (\ref{growthcondition}) implies that 
for all $j=1,\dots, n$, there exist two constants $\delta_j =\delta_j(P_j)>0$ and $C_j=C_j(P_j)>0$ such that 
\begin{equation}\label{growthconditionbis}
 P_j(x_1,\dots,x_j) \geq C_j~ (x_1+\dots+x_j)^{\delta_j} \;\;\mbox{for all}\;\; \x=(x_1,\dots,x_j)\in [1,\infty)^j.
\end{equation}
Therefore it follows that $\zeta_n(\s; \P)$ converges absolutely in the domain 
\begin{align}\label{def_DnP}
\mathcal D_n(\P):=\{\s=(s_1,\dots, s_n) \in \C^n \mid \Re (\sum_{i=j}^n \delta_i s_i) > n+1-j \quad \mbox{for all}\;\; j=1,\dots, n\}.
\end{align}
In fact, \eqref{growthconditionbis} implies 
$$\zeta_n(\s; \P)\ll 
\sum_{ m_1, \dots, m_n \geq 1}\prod_{j=1}^n(m_1+\cdots+m_j)^{-\delta_j\Re s_j},$$
and the right-hand side is convergent in the region $\mathcal D_n(\P)$
(see \cite[Theorem 3]{mat-illinois}).

Assume that for any $j=1,\dots, n$, $P_j$ satisfies the assumption 
\begin{equation*}\label{h0S}
(H_0S)\qquad 
\frac{\partial^{\alphab} P_j}{P_j} {\mbox { is bounded in }} [1,\infty)^j \quad
{\mbox{for all}} \;\alphab \in \N_0^j, 
\end{equation*}
where $\partial^{\alphab}=\partial^{\alpha_1}\cdots\partial^{\alpha_j}$
for $\alphab=(\alpha_1,\ldots,\alpha_j)$.
The method of \cite{essouabriThesis} and \cite{essouabriFourier} (see Remark 2 in page 74 of \cite{essouabriThesis}) implies that 
$\s \mapsto \zeta_n(\s; \P)$ has a meromorphic continuation to the whole space $\C^n$ and that there exists a finite set $I(\P)\subset \N_0^n$ and nonnegative integers 
$d_{\alphab}$ $(\alphab \in I(\P))$ 
such that the possibles poles are located in the set 
$$\mathcal P (\P):=\bigcup_{\alphab \in I(\P)}\bigcup_{k\in \N_0} \left\{ \s=(s_1,\dots, s_n)\in \C^n \mid \langle \s, \alphab \rangle = d_{\alphab}-k\right\}.$$

Our main result on \eqref{NLEZMzetadef} (Theorem \ref{mainsection3}) shows that if $P_j$ satisfies the above $(H_0S)$ for all $j=1,\dots, n-1$ and if $P_n$ is elliptic homogeneous then for any $\Nb=(N_1,\dots, N_n) \in \N_0^n$, 
the limit 
$$\zeta_n^{\e_n}(-\Nb; \P) :=\lim_{t\rightarrow 0} \zeta_n (-\Nb +t \e_n; \P)
=\lim_{t\rightarrow 0}\zeta_n\left((-N_1,\dots, -N_{n-1}, -N_n+t); \P\right),$$
where $\e_n=(0,\ldots,0,1)$, 
exists and can be written as a closed formula in terms of $\Nb$, the classical Bernoulli numbers and a finite number of ``periods'' (in the sense of Kontsevich-Zagier \cite{KonZag}) which depend explicitly on the polynomials $P_1,\dots, P_n$.
These periods can be interpreted as multivariate analogs of the values of the Euler gamma function at rational numbers. 

We will state Theorem \ref{mainsection3} and its corollary in Section 
\ref{sec-firstmaintheorem}.
In order to prove Theorem \ref{mainsection3}, we will first evaluate the values of
Mahler's series at non-positive integers (Theorem \ref{mainbissection3}) in
Section \ref{sec-mahlerseries}, and then prove Theorem \ref{mainsection3} in
Section \ref{sec-proof-firstmaintheorem}.

The important issue here is that {\it these periods are not necessary rational numbers and therefore (in contrast with the linear case)  the regularized values $\zeta_n^{\e_n}(-\Nb; \P)$ are not necessary in the field generated over $\Q$ by the coefficients of the polynomials $P_j$ and the direction $\e_n$.}     Some examples are given in Section \ref{exampletranscendence}.

The method developed in Section \ref{sec-mahlerseries} is influenced by the idea of 
the work of
Friedman and Pereira \cite{fredman}, so there is some common feature shared with our
previous work \cite{GMZV}.
However the method here is not a direct generalization of the method in \cite{GMZV}.
As a consequence, the formulas obtained in \cite{GMZV} and in the present paper do not
coincide, and so, comparing those two formulas we can obtain certain non-trivial
relations among Bernoulli numbers.    This point will be discussed in the last
section.

In the following sections, the empty sum is to be understood as zero.

\section{Values of $\zeta_{n, \d, \gammab}(\s)$ at $n$-tuples of non-positive integers}\label{sec-firstmainextension}

Let $n\in \N$, $\d=(d_1,\dots, d_n)\in \N^n$ and $\gammab =(\gamma_1,\dots, \gamma_n)\in \C^n$ be such that $\Re (\gamma_j)>0$ for all $j=1,\dots,n$.

In this section we state the results for the series $\zeta_{n, \d, \gammab} (\s)$
defined by \eqref{NEZMzetadef}.
The following result gives some basic properties of $\zeta_{n, \d, \gammab} (\s)$.
\begin{prop}\label{basicproperties1}
\be
\item The multiple zeta function $\zeta_{n, \d, \gammab} (\s)$ converges absolutely and uniformly  in any compact subset of the domain 
$$
\mathcal D_{n,\d}(0):=\{\s=(s_1,\dots, s_n)\in \C^n \mid \Re (s_j+\dots+ s_n) >\frac{1}{d_j}+\dots +\frac{1}{d_n} ~~(1\leq j\leq n)\},
$$
\item $\zeta_{n, \d, \gammab} (\s)$ has meromorphic continuation to the whole complex space $\C^n$ whose possible singularities are located on the union of the hyperplanes 
$$s_j+\dots+ s_n =\frac{1}{d_j}+\frac{\eps_{j+1}}{d_{j+1}}+\dots +\frac{\eps_n}{d_n} -k_j \quad (1\leq j\leq n),$$
where
$k_j\in \N_0$ and $\eps_{j+1},\dots, \eps_n \in \{0,1\}$.
\item Assume that for all $ j=1,\dots,n$ and $\eps_{j+1},\dots,\eps_n \in \{0,1\}$:
\begin{equation}\label{djsassumption}
\frac{1}{d_j} +\sum_{k=j+1}^n \frac{\eps_k}{d_k} \not \in \N.
\end{equation} 
Then, for all $\Nb=(N_1,\dots,N_n) \in \Z^n$, $\s=\Nb$ is a regular point of $\zeta_{n, \d, \gammab} (\s)$.
\ee
\end{prop}
Point 3 of Proposition \ref{basicproperties1} implies that under assumption (\ref{djsassumption})  the $n-$tuples of integers are regular point of $\zeta_{n, \d, \gammab} (\s)$.    This is an important feature different from the linear case.

Our following first main result gives a relation among those values at integer points.
\begin{theorem}\label{main-main}
Let $n\geq 2$, $\d=(d_1,\dots, d_n)\in \N^n$ and $\gammab =(\gamma_1,\dots, \gamma_n)\in \C^n$ be such that $\Re (\gamma_j)>0$ for all $j=1,\dots,n$.
Assume that the $d_j$s satisfy  the assumption (\ref{djsassumption}).
Then, for all $\Nb =(N_1,\dots, N_n)\in \Z^n$ such that $N_n\leq 0$, we have
\begin{eqnarray}\label{mainmainFormula}
\zeta_{n, \d, \gammab} (\Nb)&=&-\frac{1}{2}\zeta_{n-1, \d', \gammab'}\left(N_1,\dots,N_{n-2}, N_{n-1}+N_n\right) \nonumber\\
& -& \!\!\!\sum_{1\leq k\leq \left[(1-d_n N_n)/2\right]\atop d_n\mid 2k-1} \frac{B_{2k}}{2k}{-N_n \choose (2k-1)/d_n}\gamma_n^{(2k-1)/d_n}\nonumber\\
&&\times\zeta_{n-1, \d', \gammab'}\left(N_1,\dots,N_{n-2}, N_{n-1}+N_n +\frac{2k-1}{d_n}\right),
\end{eqnarray}
where $\d'=(d_1,\dots,d_{n-1})$ and $\gammab'=(\gamma_1,\dots,\gamma_{n-1})$.\par
\end{theorem}

Since 
$d_n | (2k-1)$, we see that
$$\left(N_1,\dots,N_{n-2}, N_{n-1}+N_n +\frac{2k-1}{d_n}\right)\in (-\N_0)^{n-1}$$ 
for $\Nb =(N_1,\dots,N_n) \in (-\N_0)^n$ and 
$k\in\{1,\dots, \left[(1-d_n N_n)/2\right]\}$.   Therefore,
Theorem \ref{main-main} gives a simple method to compute the values of $\zeta_{n, \d, \gammab} (\Nb)$ ($\Nb  \in (-\N_0)^n$) by induction on $n$.

In particular, in some special cases, Theorem \ref{main-main} gives simple closed forms
of $\zeta_{n, \d, \gammab} (\Nb)$, which we state as the following corollaries.

Let $\zeta(s)$ be the Riemann zeta-function.   It is well known that 
$\zeta(-N)\in \Q$ for all $N\in \N_0$ and 
\begin{equation} \label{TrivialRiemannZeros} 
\zeta (0)=-\frac{1}{2}\quad {\mbox { and }} \quad \zeta(-N) =0 {\mbox { for all even positive  integer }} N.
\end{equation}
With our notations here (see (\ref{NEZMzetadef})) the property (\ref{TrivialRiemannZeros}) can be reformulated as follows:
$$\zeta_{1,2,1} (0)=-\frac{1}{2}\quad {\mbox { and }} \quad \zeta_{1,2,1}(-N) =0 {\mbox { for all positive integer }} N.$$
Our first corollary of Theorem \ref{main-main} extends these properties  to the multivariable setting as follows:
\begin{coro}\label{main1}
Let $n\in \N$, $\d=(d_1,\dots, d_n)\in \N^n$ and $\gammab =(\gamma_1,\dots, \gamma_n)\in \C^n$ be such that $\Re (\gamma_j)>0$ for all $j=1,\dots,n$.
Assume that the $d_j$s satisfy  the assumption (\ref{djsassumption}).\\
Then, 
\be 
\item For all $\Nb=(N_1,\dots, N_n)\in \N_0^n$,  $\s=-\Nb$ is a regular point of $\zeta_{n, \d, \gammab}$ and 
$\zeta_{n, \d, \gammab} (-\Nb)$ lies in the field generated over $\Q$ by the coefficients $\gamma_1,\dots,\gamma_n$;
\item If $d_2,\dots,d_n$ are even integers, then for all $\Nb=(N_1,\dots, N_n)\in \N_0^n$,
 $$\zeta_{n, \d, \gammab} (-\Nb)=\left(-\frac{1}{2}\right)^{n-1} \gamma_1^{|\Nb|}(-1)^{d_1|\Nb|}\frac{B_{d_1|\Nb|+1}}{(d_1|\Nb|+1)};$$
 \item In particular, if $d_1,\dots,d_n$ are even integers, then 
$$ \zeta_{n, \d, \gammab} (\zerob_n)=\left(-\frac{1}{2}\right)^n,$$ 
where $\zerob_n=(0,\ldots,0)\in \N_0^n$, and
$$\zeta_{n, \d, \gammab} (-\Nb) =0 \;\;
{\mbox {\rm for all }}\; \Nb=(N_1,\dots, N_n)\in \N_0^n \setminus\{\zerob\}.$$
\ee
\end{coro}

If $d_2,\dots,d_n$ are not all even integers the expression of $\zeta_{n, \d, \gammab}$ at $n$-tuples of non-positive integers are more complicated. 
However, we have the following partial result:
\begin{coro}\label{main2}
Let $n\in \N$, $\d=(d_1,\dots, d_n)\in \N^n$ and $\gammab =(\gamma_1,\dots, \gamma_n)\in \C^n$ be such that $\Re (\gamma_j)>0$ for all $j=1,\dots,n$.
Assume that the $d_j$s satisfy  the assumption (\ref{djsassumption}).   Then
\be
\item
$$ \zeta_{n, \d, \gammab} (\zerob_n)=\left(-\frac{1}{2}\right)^n;$$
\item 
$$\zeta_{n, \d, \gammab}(\zerob_{n-1},-1)=\left(-\frac{1}{2}\right)^{n-1}
\left(\frac{(-1)^{d_1}B_{d_1+1}}{d_1+1}\gamma_1-\sum_{j=2}^n\frac{B_{d_j+1}}{d_j+1}
\gamma_j\right);$$
\item
\begin{align*}
&\zeta_{n, \d, \gammab}(\zerob_{n-1},-2)=\left(-\frac{1}{2}\right)^{n-1}
\frac{B_{2d_1+1}}{2d_1+1}\gamma_1^2\\
&\qquad-2\left(-\frac{1}{2}\right)^{n-2}\sum_{k=2}^n\frac{B_{d_k+1}}{d_k+1}\gamma_k
\left(\frac{(-1)^{d_1}B_{d_1+1}}{d_1+1}\gamma_1
-\sum_{j=2}^{k-1}\frac{B_{d_j+1}}{d_j+1}\gamma_j\right).
\end{align*}
\ee
\end{coro}

{\bf Remark:}
The values of  
$\zeta_{n, \d, \gammab}(\zerob_{n-2},-1,0)$ and $\zeta_{n, \d, \gammab}(\zerob_{n-2},-2,0)$
can also be computed by using points 2 and 3 and the formulas
\begin{align*}
\zeta_{n, \d, \gammab}(\zerob_{n-2},-\ell,0)&=-\frac{1}{2}\zeta_{n-1, \d', \gammab'}
(\zerob_{n-2},-\ell)
\quad (\ell=1,2)
\end{align*}
which immediately follows from \eqref{mainmainFormula}.
\medskip

The following corollary of Theorem \ref{main-main} is more intriguing (and we hope also interesting). It gives a link between the values of Riemann zeta-function at 
(odd or even) positive integers and the values of a double zeta-function at mixed pairs of integers. 
\begin{coro}\label{intriguing}
Let $d_1\in \N\setminus\{1\}$ and $d_2\in \N$ be such that $\frac{1}{d_1}+\frac{1}{2d_2}\not \in \N$.\\
Then, for all $N\in \N_0$, 
$$\zeta (d_1) = -2 ~\zeta_{2,(d_1,2 d_2),(1,1)}(1+N, -N).$$
\end{coro}

{\bf Remark:}
Since $\zeta(d_1)$ is transcendental at least when $d_1$ is even, 
Corollary \ref{intriguing} especially implies that for $\d$ satisfying the assumption (\ref{djsassumption}), {\it if the components of $\Nb\in \Z^n$ are not all non-positive, the value $\zeta_{n, \d, \gammab} (\Nb)$ is not necessary in the field generated over $\Q$ by the coefficients $\gamma_j$}.    This is contrary to the situation described
in point 1 of Corollary \ref{main1}.
\medskip

The following result shows that if the assumption (\ref{djsassumption}) does not hold, then even $n-$tuples of non-positive integers can lie on the singular locus of 
 $\zeta_{n, \d, \gammab}(\s)$. However the directional limits exist, and unlike the
 classical case, they can be highly transcendental!
 
 \begin{theorem}\label{irrationalityTh}
 Let $n\in \N$ be such that $n\geq 3$. Let $\gammab =(\gamma_1,\dots,\gamma_n)$ be such that $\Re \gamma_j >0$ for all $j=1,\dots,n$. 
 Let $\d=(d_1,\dots,d_n) \in \N^n$. Assume that for $j=1,\dots,n$ and  $\eps_{j+1},\dots,\eps_n \in \{0,1\}$:
\begin{equation}\label{ira-assumption}
 \frac{1}{d_j} +\sum_{k=j+1}^n \frac{\eps_k}{d_k}  \in \N \quad {\mbox{ if and only if }} \quad j=2 {\mbox{ and }} \eps_3 = \dots= \eps_n =1.
\end{equation} 
Denote $b:= \sum_{k=2}^n d_k^{-1}$ {\rm(}which is $\in \N$ by 
\eqref{ira-assumption}{\rm)}.
Let $\Nb=(N_1,\dots,N_n) \in \N_0^n$. Let $\thetab =(\theta_1,\dots, \theta_n)\in \C^n$ be such that $\theta_2+\dots+\theta_n\neq 0$ and $\theta_n\neq 0$. 
Then,
the directional limit $\ds  \zeta_{n, \d, \gammab}^{\thetab}(-\Nb):=\lim_{t\rightarrow 0,~ t\in \C\setminus\{0\}}  \zeta_{n, \d, \gammab}(-\Nb +t\thetab)$ exists and is given by 
\begin{align}\label{lgammavalues-zetavalues}
\zeta_{n, \d, \gammab}^{\thetab}(-\Nb)=& C(\Nb, \d)~B_{d_1(|\Nb|+b)+1}~\left(\gamma_1^{|\Nb|+b}\prod_{j=2}^n \gamma_j^{-1/d_j}\right) \left(\frac{\theta_n}{\theta_2+\dots+\theta_n}\right)~\prod_{j=2}^n \Gamma \left(\frac{1}{d_j}\right) \nonumber \\
& + H_{n,\d,\gammab}(\Nb),  
\end{align}
where
\begin{align*}
C(\Nb,\d) :=& \prod_{j=3}^n \left(\prod_{u=-N_{j-1}-\dots-N_n}^{-N_{j}-\dots-N_n-1}\left(u-\sum_{k=j}^n \frac{1}{d_k}\right)\right)\nonumber \\
& \quad \times 
\frac{(-1)^{N_2+\dots+N_{n-1}+b+d_1(|\Nb|+b)} N_n!}{(N_2+\dots+N_n+b)! (\prod_{j=2}^n d_j) \left(d_1(|\Nb|+b)+1\right)} \in \Q\setminus\{0\};
\end{align*}
\begin{align*}
&H_{n,\d,\gammab}(\Nb):=-\frac{1}{2}\zeta_{n-1, \d', \gammab'}(-N_1,\dots,-N_{n-2}, -N_{n-1}-N_n)\nonumber  \\
&\; -\sum_{k=1\atop d_n \mid 2k-1}^{\left[(1+d_nN_n)/2\right]} \frac{B_{2k}}{2k} {N_n\choose (2k-1)/d_n} \gamma_n^{(2k-1)/d_n}\nonumber\\
&\; \times  \zeta_{n-1, \d', \gammab'}\left(-N_1,\dots,-N_{n-2}, -N_{n-1}-N_n+(2k-1)/d_n\right) \in \Q(\gamma_1,\dots,\gamma_n),
  \nonumber
 \end{align*}
 where $\d':=(d_1,\dots,d_{n-1})$ and $\gammab':=(\gamma_1,\dots,\gamma_{n-1})$.
\end{theorem}
{\bf Remarks:} 
\be
\item The fact that $C(\Nb,\d)\in \Q\setminus\{0\}$ follows from assumption ({\ref{ira-assumption}}), because \eqref{ira-assumption} implies that 
$\sum_{k=j}^n d_k^{-1}\notin \N$ for $j\geq 3$. The fact that  
$H_{n,\d,\gammab}(\Nb)\in \Q(\gamma_1,\dots,\gamma_n)$ follows from point 1 of Corollary \ref{main1} since $\d'=(d_1,\dots,d_{n-1})$ satisfies  $\frac{1}{d_j} +\sum_{k=j+1}^{n-1} \frac{\eps_k}{d_k} \not \in \N$ for $ j\in\{1,\dots,n-1\}$ and $\eps_{j+1},\dots,\eps_{n-1} \in \{0,1\}$;
\item
If  $d_1(|\Nb|+b)$ is an odd integer, then $B_{d_1(|\Nb|+b)+1}\neq 0$, and hence
Theorem \ref{irrationalityTh} shows (because of the existence of the factor
$\theta_n/(\theta_2+\dots+\theta_n)$)
that $\s=-\Nb$ indeed lies on the singular locus of $\zeta_{n,\d,\gammab}(\s)$ and is a point of indeterminacy.
\ee

It is to be stressed that {\it Theorem \ref{irrationalityTh} gives a link between the diaphantine properties of the multiple zeta values of $\zeta_{n, \d, \gammab}$ at 
$n-$tuples of non positive integers and the important problem in the transcendental theory concerning the diophantine properties of the values of  Euler's gamma function at rational points.}
In fact,
$\Gamma (1/2)=\sqrt{\pi}$ is a transcendental number, and $\Gamma(1/3)$, $\Gamma(1/4)$ and $\Gamma(1/6)$ are also transcendental (Chudnovsky \cite{chu} \cite{chu2}; see also the introduction of
\cite{chu3}).    Therefore we deduce immediately from Theorem \ref{irrationalityTh} the following result.    Let $\overline{\Q}$ be the set of all algebraic numbers.

\begin{coro}\label{irrationalityCoro}
Let $q\in \{2,3,4,6\}$ and $d_1\in \N$ be such that $d_1>q$. Set $n=q+1$  and $\d=(d_1, q,\dots,q)\in \N^n$.
Let $\gammab =(\gamma_1,\dots,\gamma_n)\in {\overline{\Q}}^n$ be such that $\Re \gamma_j >0$ $\forall j=1,\dots,n$. 
Let $\thetab =(\theta_1,\dots, \theta_n)\in {\overline{\Q}}^n$ be such that $\theta_2+\dots+\theta_n\neq 0$ and $\theta_n\neq 0$.
Then, 
for all $\Nb=(N_1,\dots, N_n) \in \N_0^n$ such that $d_1(|\Nb|+b)$ is an odd integer, 
$$\zeta_{n, \d, \gammab}^{\thetab}(-\Nb)=\zeta_{n, \d, \gammab}^{\thetab}(-N_1,\dots,-N_n)$$ 
is a transcendental number.
\end{coro}

{\bf Proof:} It is enough to check that under the conditions $d_1>q$ and $n=q+1$, 
$\d=(d_1, q,\dots,q)$ satisfies \eqref{ira-assumption}, and apply Chudnovsky's result
to the gamma factors on the right-hand side of \eqref{lgammavalues-zetavalues}. \qed

\section{Values of $\zeta_n(\s; \P)$ at non-positive integers}\label{sec-firstmaintheorem}

Hereafter, for vectors $\a=(a_1,\ldots,a_n)$ and $\b=(b_1,\ldots,b_n)$ we write
$|\a|=a_1+\cdots+a_n$, $\a!=a_1!\cdots a_n!$, and
$\a^{\b}=a_1^{b_1}\cdots a_n^{b_n}$.
The inequality $\a\leq \b$ means $a_i\leq b_i$ for all $i=1,\ldots,n$.
Also we introduce the following notations.
Let $n, d\in \N$ and  $q, N \in \N_0$.
\be
\item Define for any $\betab=(\beta_1,\ldots,\beta_n) \in \N_0^n$, 
$$I_N(\betab):=\big\{\alphab=(\alpha_1,\ldots,\alpha_d) \in \N_0^d;~~ \sum_{k=1}^d k \alpha_k +|\betab| =d N +q+n\big\};$$
\item Define for any $k=1,\dots, d$,
$\Delta_k^n =\{\gammab=(\gamma_1,\ldots,\gamma_n) \in \N_0^n;~~ |\gammab|=k\},$
hence $|\Delta_k^n|={n+k-1 \choose n-1}$;
\item 
Define for any $\alphab =(\alpha_1,\ldots,\alpha_d)\in \N_0^d$,
$$V(\alphab):=\left\{\u=(\u_1,\dots, \u_d)
;~ \u_k =\left(u_{k,\gammab}\right)_{\gammab \in \Delta_k^n} \in 
\N_0^{{n+k-1 \choose n-1}},
 |\u_k| = \alpha_k ~~(1\leq k\leq d)\right\};$$
\item Define for any $\u =(\u_1,\dots, \u_d)\in \prod_{k=1}^d \N_0^{{n+k-1 \choose n-1}}$,
$$g(\u):=\left(g_1(\u),\dots,  g_n(\u)\right) {\mbox { where }}    g_i(\u)  :=\sum_{k=1}^d \sum_{\gammab \in \Delta_k^n} u_{k,\gammab} 
\gamma_i ~(1\leq i\leq n);$$
\item 
$\hat{\y}(i)=(y_1,\dots, y_{i-1}, 1, y_{i+1},\dots, y_n)$ 
for any $y_1,\dots,y_{i-1},y_{i+1},\dots, y_n$ and any $i=1,\dots,n$;
\item For any polynomial $P\in \R[X_1,\dots, X_n]$, any $\alphab \in \N_0^d$, any $\u \in V(\alphab)$ and any $i\in\{1,\dots, n\}$, we define the polynomial $P^i_{\alphab, \u}$ in 
$n-1$ variables $y_1,\dots, y_{i-1}, y_{i+1},\dots, y_n$ by 
$$P^i_{\alphab, \u}\left(\hat{\y}(i)\right):=
\frac{\alphab !}{\left(\prod_{k=1}^d \u_k !\right) 
}
\prod_{k=1}^d \prod_{\gammab \in \Delta_k^n}\left(\frac{\partial^{\gammab}P\left(\hat{\y}(i)\right)}{\gammab !}\right)^{u_{k,\gammab}}.$$
\ee 

For any elliptic and homogeneous polynomial $P\in \R[X_1,\dots, X_n]$ of degree $d$ and any polynomial $Q\in \R[X_1,\dots, X_n]$ of degree $q$, we define 
the following integrals, which are periods in the sense of Kontsevich-Zagier
\cite{KonZag} when the coefficients of relavant polynomials are rational:

\begin{definition}\label{period_def}
Let $\alphab \in \N_0^d$, $\u \in V(\alphab)$, $\betab \in \N_0^n$ and 
$i\in \{1,\dots, n\}$.
The integral (period) $K_i(P;Q; N;\alphab, \u, \betab)$ is defined by
\begin{align}\label{period_def}
K_i(P;Q; N; \alphab, \u, \betab):=\int_{[0,1]^{n-1}} 
P\left(\hat{\y}(i)\right)^{N-|\alphab|}P_{\alphab, \u}^i\left(\hat{\y}(i)\right)
(\partial^{\betab}Q\left(\hat{\y}(i)\right)) 
\prod_{k=1\atop k\neq i}^n dy_k.
\end{align}
\end{definition}
With these notations, the main result in this section is the following:
\begin{theorem}\label{mainsection3}
Consider a polynomial $ P_j\in \R[X_1,\dots, X_j]$ in $j$ variables for any $j=1,\dots, n$.
Assume that the assumptions (\ref{growthcondition}), (\ref{forsimplicity}) hold. 
Assume also that for all $j=1,\dots n-1$ the polynomial $P_j$ satisfies the assumption $(H_0S)$ and that the polynomial 
$P_n$ is elliptic and homogeneous of degree $d\geq 1$.
Denote by $\zeta_n(\s; \P)$ the meromorphic continuation of 
$$\s=(s_1,\dots, s_n) \mapsto \sum_{ m_1, \dots, m_n \geq 1} \frac{1}{\prod_{j=1}^n P_j(m_1, \dots, m_j)^{s_j}}$$
to the whole complex space $\C^n$.
Then,
 for any $\Nb=(N_1,\dots, N_n) \in \N_0^n$, the limit 
 $$\zeta_n^{\e_n}(-\Nb; \P) :=\lim_{t\in\C,t\rightarrow 0} \zeta_n (-\Nb +t \e_n; \P)= \lim_{t\in\C,t\rightarrow 0} \zeta_n\left((-N_1,\dots, -N_{n-1}, -N_n+t); \P\right)$$
 exists and is given  by 
 \begin{eqnarray}\label{explicitvalue}
 \zeta_n^{\e_n}(-\Nb; \P) &=&
 \sum_{\betab \in \N_0^n\atop |\betab|\leq q_{\Nb}} \sum_{\alphab \in \N_0^d \atop \sum_{k=1}^d k\alpha_k +|\betab|=q_{\Nb}+n} \sum_{\u \in V(\alphab)}
 \frac{(-1)^{|\alphab|} (|\alphab|-1)!}{d ~\alphab ! ~\betab !} \nonumber\\
 &&\times\left(\sum_{i=1}^n K_i(P_n;Q_{\Nb}; 0; \alphab, \u, \betab)\right) 
 \prod_{i=1}^n {\widetilde{B}}_{g_i(\u)+\beta_i},
 \end{eqnarray}
 where 
 \be
 \item ${\widetilde{B}}_k :=B_k$ (the classical Bernoulli number) for all $k\neq 1$ and 
 ${\widetilde{B}}_1:=-B_1=\frac{1}{2}$,
 \item $Q_{\Nb} = \prod_{j=1}^n P_j(X_1,\dots, X_j)^{N_j}$ and $q_{\Nb}=\deg Q_\Nb =\sum_{j=1}^n N_j \deg P_j$.
 \ee
\end{theorem}

{\bf Remark}$\;$ Recently some authors (such as \cite{AIK})
prefer to define Bernoulli numbers in a slightly
different way from \eqref{def_Bernoulli}, that is, by
$$
\frac{Xe^X}{e^X-1}=\sum_{n=0}^{\infty}\frac{\widetilde{B}_n}{n!}X^n.
$$
This $\widetilde{B}_n$ is exactly the same as $\widetilde{B}_n$ above.
\bigskip

An interesting illustration of Theorem \ref{mainsection3} is given by its following corollary.
Define the generalized gamma factor $G_{n-1}(m; \mub)$ by 
$$G_{n-1}(m; \mub):=\int_{(0,1)^{n-1}} \frac{\prod_{i=1}^{n-1} t_i^{\mu_i -1}~dt_1\dots dt_{n-1}}{\left(1+t_1+\dots+t_{n-1}\right)^m}$$
for any $m\in \N_0$ and any $\mub=(\mu_1,\ldots,\mu_{n-1}) \in (0,\infty)^{n-1}$. 
\begin{coro}\label{mainsection3coro}
Consider for any $j=1,\dots, n$ a polynomial $ P_j\in \R[X_1,\dots, X_j]$ in $j$ variables.
Assume that the assumptions (\ref{growthcondition}), (\ref{forsimplicity}) hold. 
Assume also that for all $j=1,\dots n-1$ the polynomial $P_j$ satisfies the assumption $(H_0S)$ and that 
$$P_n=X_1^d+\dots+X_n^d \quad {\mbox { where }} d\geq 1.$$
Then,
 \begin{eqnarray*}
 \zeta_n^{\e_n}(-\Nb; \P) &=&
 \sum_{\betab \in \N_0^n\atop |\betab|\leq q_{\Nb}} \sum_{\nub \in \N_0^n\atop \nub\leq \betab} \sum_{\alphab \in \N_0^d \atop \sum_{k=1}^d k\alpha_k =|\nub|+n} 
 \sum_{\gammab=(\gammab^1,\dots, \gammab^d) \in (\N_0^n)^d \atop |\gammab^k|=\alpha_k ~(1\leq k\leq d)}
 \frac{(-1)^{|\alphab|} (|\alphab|-1)! {\betab \choose \nub} \partial^{\betab}Q_{\Nb}(\zerob)}{d^n ~\betab !~\prod_{k=1}^d \gammab^k !} \\
 && \left[\sum_{i=1}^n G_{n-1}\left(|\alphab|; \mub^i (\nub; \gammab)\right)\right]
 \times \prod_{i=1}^n {\widetilde{B}}_{\beta_i-\nu_i +\sum_{k=1}^n k \gamma^k_i},
 \end{eqnarray*}
 where 
for any $\nub=(\nu_1,\ldots,\nu_n) \in \N_0^n$ and any $\gammab=(\gammab^1,\dots, \gammab^d) \in (\N_0^n)^d$ with $\gammab^k=(\gamma_1^k,\ldots,\gamma_n^k)\in\N_0^n$,
 $$\mub^i (\nub; \gammab)=\left(\mu^i_1(\nub; \gammab),\dots,\mu^i_{i-1}(\nub; \gammab),\mu^i_{i+1}(\nub; \gammab),\dots, \mu^i_{n}(\nub; \gammab)\right)$$
with
 $\mu^i_j(\nub; \gammab)=\frac{1}{d}\left(1+\nu_j+\sum_{k=1}^d (d-k) \gamma^k_j\right).$
\end{coro}

\section{Proof of Proposition \ref{basicproperties1}}\label{sec-proofstheorems123}
\subsection{Three elementary lemmas}
The following three lemmas are elementary but useful for our proofs.
\begin{lemma}\label{integral}
Let $d>0$. Let $a, b \in \C$ be such that $\Re a >0$ and $\Re b >0$. Let $s\in \C$ be such that $\Re(s)>1/d$. Then 
$$\int_0^\infty (b+a x^d)^{-s} ~dx = \frac{1}{d ~a^{1/d} ~b^{s-1/d}}\frac{\Gamma(s-1/d) \Gamma(1/d)}{\Gamma(s)}.$$
\end{lemma}
{\bf Proof of Lemma \ref{integral}:} 
Assume first that $a, b\in (0,\infty)$. By using the change of variables 
$x=(b/a)^{1/d} (t^{-1}-1)^{1/d}$, we obtain that\\
\begin{eqnarray*}
\int_0^\infty (b+ax^d)^{-s} ~dx &=&\frac{1}{d b^{s-1/d} a^{1/d}}\int_0^1 (1-t)^{1/d -1} t^{s-1/d -1}~dt \\
&=& \frac{ B(1/d, s-1/d)}{d b^{s-1/d} a^{1/d}}=
\frac{\Gamma(s-1/d) \Gamma(1/d)}{d b^{s-1/d} a^{1/d} \Gamma(s)}.
\end{eqnarray*}
By using in addition the fact that the function $(a,b)\mapsto \int_0^\infty (b+a x^d)^{-s} ~dx$ is holomorphic in the domain $\{(a,b)\in \C^2\mid \Re a>0,~\Re b>0\}$, we conclude by analytic continuation 
that the lemma holds for all $a, b\in\C$ such that $\Re a >0$ and $\Re b >0$. \qed

\begin{lemma}\label{EulerMaclaurin}
Let $K\in \N$ and $f\in \mathcal C^{(2K)}\left([0,\infty),\C\right)$. 
Assume that $\int_0^\infty |f(x)|~dx <+\infty$,  $\int_0^\infty |f^{(2K)}(x)|~dx <+\infty$ and $\lim_{x\rightarrow \infty} f^{(k)}(x)=0$ for all $k=0,\dots,2K-1$.
Then, the series $\sum_{m\geq 1} f(m)$ is convergent and 
$$\sum_{m=1}^\infty  f(m)=\int_0^\infty f(x)~dx -\frac{f(0)}{2}-\sum_{k=1}^K \frac{B_{2k}}{(2k)!}  f^{(2k-1)}(0)-\frac{1}{(2K)!}\int_0^\infty f^{(2K)}(x) B_{2K}(\{x\})~dx,$$
where $(B_k)_{k\geq 0}$ is the sequence of Bernoulli numbers and $(B_k(\cdot))_{k\geq 0}$ is the sequence of Bernoulli polynomials with $\{x\}$ the fractional part of $x$.
\end{lemma}
{\bf Proof of Lemma \ref{EulerMaclaurin}:} it follows easily from the classical Euler-Maclaurin formula. \qed

\begin{lemma}\label{faderivatives}
Let $d\in \N$ and $\s \in \C$. Let $a, b\in \C$ be such that $\Re a >0$ and $\Re b >0$. Set $\delta:=(\Re b/\Re a)^{1/d}$. Define the function 
$f_{a,b,d,s}:(-\delta,\infty)\rightarrow \C$ by
 $f_{a,b,d,s}(x)=(b+ a x^d)^{-s}$ for any $x>-\delta$.
Then, 

1. $f_{a,b,d,s}$ is $\mathcal C^\infty$ on $(-\delta, \infty)$.

2. For $k\in \N_0$ and $x>-\delta$,
$$f_{a,b,d,s}^{(k)}(x)=k!\sum_{\alphab \in \N_0^d\atop \sum_{j=1}^d j\alpha_j =k}{-s\choose |\alphab|}\frac{|\alphab| !}{\alphab !}\left(\prod_{j=1}^d {d\choose j}^{\alpha_j}\right)
a^{|\alphab|}~x^{\sum_{j=1}^d (d-j)\alpha_j} ~(b+a x^d)^{-s-|\alphab|}.$$

3. 
In particular, $ f_{a,b,d,s}^{(k)}(0)=k! {-s\choose k/d} a^{k/d}~b^{-s-k/d}$ if $d\mid k$ and $f_{a,b,d,s}^{(k)}(0)=0$ if $d\nmid k$.
\end{lemma}
{\bf Proof of Lemma \ref{faderivatives}:}
When $x>-\delta$, we have
$\Re(b+ax^d)>0$.    Therefore 
$(b+ax^d)^{-s}=\exp(-s\log(b+ax^d))$ is well-defined and $f_{a,b,d,s}$ is 
$\mathcal C^\infty$, hence point 1.

For point 2, fix $x>-\delta$ and choose $\eps=\eps_x>0$ small enough such that  
$x+u>-\delta$ and $|a \sum_{j=1}^d {d\choose j} x^{d-j} u^j|\leq (1/2) |b+a x^d|$  for all $u\in (-\eps, \eps)$.
It follows then from Newton's binomial formula that for $u\in (-\eps, \eps)$:
\begin{eqnarray*}
f_{a,b,d,s}(x+u)&=& \left(b+a(x+u)^d\right)^{-s}=(b+ax^d)^{-s}\left(1+\frac{a\sum_{j=1}^d {d\choose j} x^{d-j} u^j}{b+ax^d}\right)^{-s}\\
&=& (b+ax^d)^{-s}\sum_{\ell=0}^\infty {-s\choose \ell} \frac{\left(\sum_{j=1}^d {d\choose j} x^{d-j} u^j\right)^\ell}{(b+a x^d)^\ell} a^{\ell}\\
&=&\sum_{\ell=0}^\infty \frac{{-s\choose \ell} a^{\ell}}{(b+a x^d)^{s+\ell}}\sum_{\alphab \in \N_0^d,~|\alphab|=\ell}\frac{\ell !}{\alphab !} 
\prod_{j=1}^d \left({d\choose j} x^{d-j} u^j\right)^{\alpha_j}\\
&=& \sum_{\alphab \in \N_0^d}\frac{{-s\choose |\alphab|} a^{|\alphab|}}{(b+a x^d)^{s+|\alphab|}}\frac{|\alphab| !}{\alphab !} 
\left(\prod_{j=1}^d {d\choose j}^{\alpha_j} \right) x^{\sum_{j=1}^d (d-j)\alpha_j } ~ u^{\sum_{j=1}^d j\alpha_j}\\
&=& \sum_{k=0}^\infty \left(\sum_{\alphab \in \N_0^d,~\sum_{j=1}^d j\alpha_j=k}\frac{{-s\choose |\alphab|} a^{|\alphab|}}{(b+a x^d)^{s+|\alphab|}}\frac{|\alphab| !}{\alphab !} 
\left(\prod_{j=1}^d {d\choose j}^{\alpha_j} \right) x^{\sum_{j=1}^d (d-j)\alpha_j }\right) u^k.
\end{eqnarray*}
Point 2 then follows from the uniqueness of Taylor's expansion of $f_{a,b}$ in a neighborhood of $x$.

Lastly, putting $x=0$, the only non-vanishing terms are under the condition
$\sum_{j=1}^d (d-j)\alpha_j=0$, which is equivalent to 
$\alpha_1=\cdots=\alpha_{d-1}=0$.  Then $d\alpha_d=k$, which implies point 3.
\qed
\subsection{Domain of absolute convergence of $\zeta_{n, \d, \gammab} (\s)$: proof of point 1 of Proposition \ref{basicproperties1}}

We will prove point 1 of Proposition \ref{basicproperties1}  by induction on $n$.
If $n=1$,  $\zeta_{1, d_1, \gamma_1} (s)=\gamma_1^{-s} \zeta (d_1 s)$ and  point 1 of Proposition \ref{basicproperties1} clearly holds. 

Assume now that $n\geq 2$ and that point 1 of Proposition \ref{basicproperties1} holds for $n-1$.
Let $\mathcal K$ be a compact subset of $\mathcal D_{n,\d}(0)$. 
Because of the definition of $\mathcal D_{n,\d}(0)$ we see that $\Re s_n >1/d_n$ for all $s_n \in \mathcal K$.   Therefore uniformly in 
$\s\in \mathcal K$ and $m_1,\dots, m_{n-1}\in \N$ we have
\begin{eqnarray*}
\sum_{m_n=1}^\infty|(\gamma_1 m_1^{d_1}+\dots+ \gamma_n m_n^{d_n} )^{-s_n}|&\ll_{\mathcal K, \gammab}&
\sum_{m_n=1}^\infty (m_1^{d_1}+\dots+m_n^{d_n})^{-\Re s_n} \\
&\leq &
 \int_{0}^{\infty} (m_1^{d_1}+\dots+m_{n-1}^{d_{n-1}}+x^{d_n})^{-\Re s_n} ~dx\\
&\ll_{\mathcal K, \d}& (m_1^{d_1}+\dots+m_{n-1}^{d_{n-1}})^{-\Re s_n+1/d_n},
\end{eqnarray*}
where the last inequality is by
Lemma \ref{integral}.
We deduce that, uniformly in $\s\in \mathcal K$,
\begin{align*}
\lefteqn{\sum_{m_1,\dots,m_n \geq 1} |\prod_{j=1}^n (\gamma_1 m_1^{d_1}+\dots+ \gamma_j m_j^{d_j} )^{-s_j}|}\\
&\ll_{\mathcal K, \d, \gammab}  \sum_{m_1,\dots,m_{n-1} \geq 1} \prod_{j=1}^{n-1} (m_1^{d_1}+\dots+  m_j^{d_j} )^{-\Re s_j}(m_1^{d_1}+\dots+m_{n-1}^{d_{n-1}})^{-\Re s_n+1/d_n}\\
&\ll_{\mathcal K, \d, \gammab} \zeta_{n-1, \d', \unb} (\Re s_1,\dots,\Re s_{n-2}, \Re s_{n-1}+\Re s_n-1/d_n)
\end{align*}
and we conclude the assertion by induction hypothesis.
This ends the proof of point 1 of Proposition \ref{basicproperties1}. \qed
\subsection{A key proposition}

Let $n\in \N$, $\d=(d_1,\dots, d_n)\in \N^n$ and $\gammab =(\gamma_1,\dots, \gamma_n)\in \C^n$ be such that $\Re (\gamma_j)>0$ for all $j=1,\dots,n$.
Define for $\delta \in \R$ the set 
\begin{align}\label{domainprol}
\lefteqn{\mathcal D_{n,\d}(\delta)}\notag\\
& :=\{\s=(s_1,\dots, s_n)\in \C^n \mid \Re (s_j+\dots+ s_n) >\left(\frac{1}{d_j}+\dots+\frac{1}{d_n}\right)+\frac{\delta}{d_n}~~(j=1,\dots, n)\}.
\end{align}
The following proposition is a key ingredient in this section:

\begin{prop}\label{mainproposition1}
 Let $K\in \N$. There exists a function $\s\mapsto \mathcal R_{n,\d, K}(\s)$ holomorphic in $\mathcal D_{n,\d}(-2K)$ such that 
for all $\s \in \mathcal D_{n,\d}(0)$, we have
\begin{eqnarray}\label{mainformula1}
\zeta_{n, \d, \gammab} (\s)&=&\frac{\Gamma(s_n-1/d_n) \Gamma(1/d_n) \gamma_n^{-1/d_n}}{d_n \Gamma (s_n)}\zeta_{n-1, \d', \gammab'}\left(s_1,\dots,s_{n-2}, s_{n-1}+s_n -1/d_n\right)\nonumber \\
& & -\frac{1}{2}\zeta_{n-1, \d', \gammab'}\left(s_1,\dots,s_{n-2}, s_{n-1}+s_n\right) \nonumber\\
& & -\sum_{k=1\atop d_n\mid 2k-1}^K \frac{B_{2k}}{2k}{-s_n \choose (2k-1)/d_n}\gamma_n^{(2k-1)/d_n}\nonumber\\
&&\qquad\times\zeta_{n-1, \d', \gammab'}\left(s_1,\dots,s_{n-2}, s_{n-1}+s_n +(2k-1)/d_n\right)\nonumber \\
& & +  \mathcal R_{n,\d, K}(\s),
\end{eqnarray}
where $\d'=(d_1,\dots,d_{n-1})$ and $\gammab'=(\gamma_1,\dots,\gamma_{n-1})$.
Furthermore, if  $\Nb =(N_1,\dots, N_n) \in \Z^n$ be such that $N_n \leq 0$ and 
\begin{align}\label{KdddN}
2K >d_n(\frac{1}{d_j}+\dots+\frac{1}{d_n}) - d_n (N_j+\dots+ N_n )
\end{align}
for all $j=1,\dots,n$, then 
$\Nb \in \mathcal D_{n,\d}(-2K)$ and $\mathcal R_{n,\d, K}(\Nb) =0$.
\end{prop}

{\bf Remark:}
Once the meromorphic continuation of $\zeta_{n, \d, \gammab} (\s)$ (point 2 of
Proposition \ref{basicproperties1} is proved, we can claim that \eqref{mainformula1} is
valid, as an identity of meromorphic functions, in the wider region
$\mathcal D_{n,\d}(-2K)$.
\medskip

{\bf Proof of Proposition \ref{mainproposition1}:}
Let $K\in \N$, $\s\in  \mathcal D_{n,\d}(0)$, and let $\m'=(m_1,\dots ,m_{n-1})\in \N^{n-1}$. Set $b(\m')=\gamma_1 m_1^{d_1} +\dots +\gamma_{n-1} m_{n-1}^{d_{n-1}}$. 
By using notations of Lemma \ref{faderivatives}, it follows from Lemma \ref{EulerMaclaurin} that 
\begin{eqnarray*}
&& \sum_{m_n =1}^\infty (\gamma_1 m_1^{d_1} +\dots +\gamma_{n} m_{n}^{d_{n}})^{-s_n}\\
&=& \sum_{m_n =1}^\infty (b(\m')+ \gamma_n m_n^{d_n})^{-s_n} =
\sum_{m_n =1}^\infty f_{\gamma_n,b(\m') ,d_n,s_n} (m_n)\\
&=& \int_0^\infty f_{\gamma_n, b(\m'),d_n,s_n}(x)~dx -\frac{f_{\gamma_n, b(\m'),d_n,s_n}(0)}{2}-\sum_{k=1}^K \frac{B_{2k}}{(2k)!}  f_{\gamma_n, b(\m'),d_n,s_n}^{(2k-1)}(0)\\
& & -\frac{1}{(2K)!}\int_0^\infty f_{\gamma_n, b(\m'),d_n,s_n}^{(2K)}(x) B_{2K}(\{x\})~dx.
\end{eqnarray*}
On the right-hand side, apply Lemma \ref{integral} to the first term, point 3 of 
Lemma \ref{faderivatives} to the third term, and point 2 of Lemma \ref{faderivatives}
to the fourth term to obtain
\begin{eqnarray}\label{preparationformula}
&& \sum_{m_n =1}^\infty (\gamma_1 m_1^{d_1} +\dots +\gamma_{n} m_{n}^{d_{n}})^{-s_n}
\notag\\
&=& \frac{1}{d_n ~\gamma_n^{1/d_n} ~b(\m')^{s_n-1/d_n}}\frac{\Gamma(s_n-1/d_n) \Gamma(1/d_n)}{\Gamma(s_n)}
- \frac{1}{2~b(\m')^{s_n}}\notag\\
& & -\sum_{k=1\atop d_n\mid 2k-1}^K \frac{B_{2k}}{2k}{-s_n \choose (2k-1)/d_n}\gamma_n^{(2k-1)/d_n} ~b(\m')^{-s_n-(2k-1)/d_n}\notag\\
& & -\sum_{\alphab \in \N_0^{d_n}\atop \sum_{j=1}^{d_n} j\alpha_j =2K}{-s_n\choose |\alphab|} G_{n,\d, K}(\m'; s_n; \alphab),
\end{eqnarray}
where
\begin{align*}
\lefteqn{G_{n,\d, K}(\m'; s_n; \alphab):=\frac{|\alphab| !}{\alphab !}\left(\prod_{j=1}^{d_n} {d_n\choose j}^{\alpha_j}\right)
\gamma_n^{|\alphab|}} \\
&\times\int_0^\infty B_{2K}(x)
~x^{\sum_{j=1}^{d_n} (d_n-j)\alpha_j} ~(b(\m')+\gamma_n x^{d_n})^{-s_n-|\alphab|}~dx.
\end{align*}
Since $\s\in \mathcal D_{n,\d}(0)$, as we already proved,
$\zeta_{n, \d, \gammab} (\s)$ converges absolutely.
Substituting \eqref{preparationformula} to the series expression of
$\zeta_{n, \d, \gammab} (\s)$, we obtain
\begin{eqnarray}\label{mainformula1bis}
\zeta_{n, \d, \gammab} (\s)&=&\frac{\Gamma(s_n-1/d_n) \Gamma(1/d_n) \gamma_n^{-1/d_n}}{d_n \Gamma (s_n)}\zeta_{n-1, \d', \gammab'}\left(s_1,\dots,s_{n-2}, s_{n-1}+s_n -1/d_n\right)\nonumber \\
& & -\frac{1}{2}\zeta_{n-1, \d', \gammab'}\left(s_1,\dots,s_{n-2}, s_{n-1}+s_n\right) \nonumber\\
& & -\sum_{k=1\atop d_n\mid 2k-1}^K \frac{B_{2k}}{2k}{-s_n \choose (2k-1)/d_n}\gamma_n^{(2k-1)/d_n}\nonumber\\
& &\qquad\times\zeta_{n-1, \d', \gammab'}\left(s_1,\dots,s_{n-2}, s_{n-1}+s_n +(2k-1)/d_n\right)\nonumber \\
& & +  \mathcal R_{n,\d, K}(\s),
\end{eqnarray}
where 
\begin{align}\label{RndKs}
\mathcal R_{n,\d, K}(\s)=-\sum_{\alphab \in \N_0^{d_n}\atop \sum_{j=1}^{d_n} j\alpha_j =2K}{-s_n\choose |\alphab|} Z_{n,\d,K, \alphab}(\s)
\end{align}
with 
$$ Z_{n,\d,K, \alphab}(\s):= 
\sum_{\m'\in \N^{n-1}}\frac{G_{n,\d, K}(\m'; s_n; \alphab) }{\prod_{j=1}^{n-1} (\gamma_1 m_1^{d_1}+\dots+ \gamma_j m_j^{d_j} )^{s_j}}.$$

Now we will prove that $\s\mapsto \mathcal R_{n,\d, K}(\s)$ is holomorphic in $\mathcal D_{n,\d}(-2K)$.
Let $\mathcal H$ be a compact subset of $D_{n,\d}(-2K)$. In particular, for all $\s\in \mathcal H$, $\Re s_n  +(2K)/d_n>1/d_n$. 
Lemma \ref{integral} implies then that we have uniformly in $\s\in \mathcal H$, $\m'\in \N^{n-1}$ and $\alphab \in \N_0^{d_n}$ such that $\sum_{j=1}^{d_n} j\alpha_j =2K$,
\begin{eqnarray*}
|G_{n,\d, K}(\m'; s_n; \alphab)| &\ll_{\mathcal H, \d, n,\gammab, K}& \int_0^\infty 
x^{ d_n |\alphab| -2K} ~(|b(\m')|+ x^{d_n})^{-\Re s_n-|\alphab|}~dx\\
& \ll_{\mathcal H, \d, n,\gammab, K}& \int_0^\infty (|b(\m')|+ x^{d_n})
^{-\Re s_n-(2K)/d_n}~dx\\
&\ll_{\mathcal H, \d, n,\gammab, K}& |b(\m')|^{-\Re s_n-(2K-1)/d_n}\\
&\ll_{\mathcal H, \d, n,\gammab, K}& (m_1^{d_1} +\dots + m_{n-1}^{d_{n-1}})^{-\Re s_n-(2K-1)/d_n}.
\end{eqnarray*}
We deduce that we have uniformly in $\s\in \mathcal H$ and $\alphab \in \N_0^{d_n}$ such that $\sum_{j=1}^{d_n} j\alpha_j =2K$,
\begin{eqnarray*}
\lefteqn{\sum_{\m'\in \N^{n-1}}\left|\frac{G_{n,\d, K}(\m'; s_n; \alphab) }{\prod_{j=1}^{n-1} (\gamma_1 m_1^{d_1}+\dots+ \gamma_j m_j^{d_j} )^{s_j}}\right| }\\
&\ll_{\mathcal H, \d, n,\gammab, K}& \sum_{\m'\in \N^{n-1}}\left(\prod_{j=1}^{n-2} ( m_1^{d_1}+\dots+ m_j^{d_j} )^{-Re s_j}\right)\\
&&\times (m_1^{d_1} +\dots + m_{n-1}^{d_{n-1}})^{-\Re s_{n-1}-\Re s_n-(2K-1)/d_n}.
\end{eqnarray*}
The set 
$$\{\left(\Re s_1,\dots,\Re s_{n-2}, \Re s_{n-1}+\Re s_n+(2K-1)/d_n\right)
 \mid \s \in \mathcal H\}$$ 
is a compact subset of $D_{n-1,\d'}(0)$, because for $s\in\mathcal H$
$$
\Re(s_j+\cdots+s_n)+\frac{2K-1}{d_n}>\frac{1}{d_j}+\cdots+\frac{1}{d_n}
-\frac{2K}{d_n}+\frac{2K-1}{d_n}=\frac{1}{d_j}+\cdots+\frac{1}{d_{n-1}}.
$$
We deduce then from point 1 of Proposition \ref{basicproperties1} that for all  $\alphab \in \N_0^{d_n}$ such that $\sum_{j=1}^{d_n} j\alpha_j =2K$, 
$\s\mapsto Z_{n,\d,K, \alphab}(\s)$ is holomorphic in $D_{n,\d}(-2K)$.
Therefore $\s\mapsto \mathcal R_{n,\d, K}(\s)$ is also holomorphic in $D_{n,\d}(-2K)$.
This implies the first assertion of Proposition \ref{mainproposition1}.

Let  $\Nb =(N_1,\dots, N_n) \in \Z^n$ be such that $N_n \leq 0$ with
\eqref{KdddN} for all $j$.
It is clear that $\Nb \in D_{n,\d}(-2K)$.
Moreover, for $\alphab \in \N_0^{d_n}$ such that $\sum_{j=1}^{d_n} j\alpha_j =2K$,
we have 
$d_n |\alphab|\geq \sum_{j=1}^{d_n} j\alpha_j =2K > -d_n N_n +1$ 
(where the last inequality is the case \textcolor{red}{$j=n$} of \eqref{KdddN}),
and therefore $|\alphab| >-N_n$.
It follows that ${-s_n\choose |\alphab|}|_{s_n=N_n}=0$. 
We conclude then from (\ref{RndKs}) that $\mathcal R_{n,\d, K}(\Nb)=0$.
This ends the proof of Proposition \ref{mainproposition1}. \qed

\subsection{Proof of points 2 and 3 of Proposition \ref{basicproperties1}}

We will first prove point 2 of Proposition \ref{basicproperties1} by induction on $n$:
If $n=1$,  $\zeta_{1, d_1, \gamma_1} (s)=\gamma_1^{-s} \zeta (d_1 s)$ and  point 2 of Proposition \ref{basicproperties1} clearly holds. 

Let $n\geq 2$. Assume that  point 2 of Proposition \ref{basicproperties1} holds for $n-1$.
By letting $K$ to infinity in (\ref{mainformula1}), we deduce then that 
$\s\mapsto \zeta_{n, \d, \gammab} (\s)$ has meromorphic continuation to $\C^n$ and that the possibles singularities are located in the union of the following  hyperplanes:
\be
\item $s_n-\frac{1}{d_n} =-k_n$ $(k_n \in \N_0)$;
\item $s_j+\dots+ s_n -\frac{1}{d_n}=\frac{1}{d_j}+\frac{\eps_{j+1}}{d_{j+1}}+\dots +\frac{\eps_{n-1}}{d_{n-1}} -k_j \quad (1\leq j\leq n-1)$, where \\
$k_j\in \N_0$ and $\eps_{j+1},\dots, \eps_{n-1} \in \{0,1\}$;
\item 
$s_j+\dots+ s_n =\frac{1}{d_j}+\frac{\eps_{j+1}}{d_{j+1}}+\dots +\frac{\eps_{n-1}}{d_{n-1}} -k_j \quad (1\leq j\leq n-1)$, where \\
$k_j\in \N_0$ and $\eps_{j+1},\dots, \eps_{n-1} \in \{0,1\}$;
\item
$s_j+\dots+ s_n +\frac{2k-1}{d_n}=\frac{1}{d_j}+\frac{\eps_{j+1}}{d_{j+1}}+\dots +\frac{\eps_{n-1}}{d_{n-1}} -k_j \quad (1\leq j\leq n-1)$, where \\
$k_j\in \N_0$, $k\in \N$ be such that $d_n \mid 2k-1$ and $\eps_{j+1},\dots, \eps_{n-1} \in \{0,1\}$.
\ee
It follows that 
$\zeta_{n, \d, \gammab} (\s)$ has meromorphic continuation to the whole complex space $\C^n$ whose possibles singularities are located in the union of the hyperplanes 
$$s_j+\dots+ s_n =\frac{1}{d_j}+\frac{\eps_{j+1}}{d_{j+1}}+\dots +\frac{\eps_n}{d_n} -k_j \quad (1\leq j\leq n),$$
where
$k_j\in \N_0$ and $\eps_{j+1},\dots, \eps_n \in \{0,1\}$.
This ends the induction argument and also the proof of point 2 of  Proposition \ref{basicproperties1}.

Now we will prove point 3  of Proposition \ref{basicproperties1}. Let $\Nb=(N_1,\dots,N_n)\in \Z^n$. Assume that $\Nb$ is a singular point of $\zeta_{n, \d, \gammab} (\s)$. It follows then from  point 2 of Proposition \ref{basicproperties1} that there exist $j\in \{1,\dots, n\}$, $\eps_{j+1},\dots, \eps_n \in \{0,1\}$ and $k_j\in \N_0$ such that 
$$N_j+\dots+ N_n =\frac{1}{d_j}+\frac{\eps_{j+1}}{d_{j+1}}+\dots +\frac{\eps_n}{d_n} -k_j.$$
Therefore $\frac{1}{d_j}+\frac{\eps_{j+1}}{d_{j+1}}+\dots +\frac{\eps_n}{d_n}\in \Z$,
but this is positive.
It follows then that 
$\frac{1}{d_j}+\frac{\eps_{j+1}}{d_{j+1}}+\dots +\frac{\eps_n}{d_n}\in \N$, 
which contradicts the assumption (\ref{djsassumption}). 
This proves point 3 and ends the proof of Proposition \ref{basicproperties1}.\qed

\section{Proof of Theorem \ref{main-main} and of Corollaries \ref{main1}, \ref{main2} and \ref{intriguing}}\label{sec-proof_th1_cor123}

{\bf Proof of Theorem \ref{main-main}:}
Let $n\in \N$, $\d=(d_1,\dots, d_n)\in \N^n$ and $\gammab =(\gamma_1,\dots, \gamma_n)\in \C^n$ be such that $\Re (\gamma_j)>0$ for all $j=1,\dots,n$.
Assume that the $d_j$s satisfy  the assumption (\ref{djsassumption}).\par
Let $\Nb =(N_1,\dots, N_n)\in \Z^n$ such that $N_n\leq 0$. Let $K\in \N$ be such that $2K >d_n(\frac{1}{d_j}+\dots+\frac{1}{d_n}) - d_n (N_j+\dots+ N_n )$ for all $j=1,\dots,n$. It follows then from 
Proposition \ref{mainproposition1} that 
\eqref{mainformula1} holds for all $\s\in D_{n,\d}(-2K)$,
moreover $\mathcal R_{n,\d, K}(\Nb) =0$ because $\s=\Nb \in \mathcal D_{n,\d}(-2K)$.

We also know from point 3 of  Proposition \ref{basicproperties1}, that $\s=\Nb$ is a regular point of $\zeta_{n, \d, \gammab} (\s)$.
Since \eqref{djsassumption} for $n-1$ is just the case $\varepsilon_n=0$ in \eqref{djsassumption} for $n$, point 3 also implies that
$\s=\Nb$ is a regular point of
$\zeta_{n-1, \d', \gammab'}\left(s_1,\dots,s_{n-2}, s_{n-1}+s_n\right)$ and \\
$\zeta_{n-1, \d', \gammab'}\left(s_1,\dots,s_{n-2}, s_{n-1}+s_n +(2k-1)/d_n\right)$ ($k\in\{1,\dots,K\}$ such that $d_n \mid 2k-1$).
Furthermore, $\s=\Nb$ is also a regular point of  $\zeta_{n-1, \d', \gammab'}\left(s_1,\dots,s_{n-2}, s_{n-1}+s_n -1/d_n\right)$.
In fact, 
if $\s=\Nb$ is a singular point of $\zeta_{n-1, \d', \gammab'}\left(s_1,\dots,s_{n-2}, s_{n-1}+s_n -1/d_n\right)$, then the point 2 of Proposition \ref{basicproperties1} implies that there exists $j\in \{1,\dots,n-1\}$, $(\eps_{j+1},\dots,\eps_{n-1})\in \{0,1\}^{n-1-j}$ and $k_j \in \N_0$ such that 
$$N_j+\dots +N_n-\frac{1}{d_n} =\left(\frac{1}{d_j}+\frac{\eps_{j+1}}{d_{j+1}}+\dots+\frac{\eps_{n-1}}{d_{n-1}}\right)-k_j.$$
This implies that $\frac{1}{d_j}+\frac{\eps_{j+1}}{d_{j+1}}+\dots+\frac{\eps_{n-1}}{d_{n-1}}+\frac{1}{d_n} \in \Z$, but this is positive, hence $\in \N$. Hence a contradiction with the assumption (\ref{djsassumption}). 

By using in addition the fact that $\frac{1}{\Gamma(s_n)}|_{s_n= N_n}=0$, we deduce then from (\ref{mainformula1}) that 
\begin{eqnarray*}
\zeta_{n, \d, \gammab} (\Nb)&=&-\frac{1}{2}\zeta_{n-1, \d', \gammab'}\left(N_1,\dots,N_{n-2}, N_{n-1}+N_n\right) \nonumber\\
& -& \!\!\!\sum_{k=1\atop d_n\mid 2k-1}^{K}\frac{B_{2k}}{2k}
\binom{-N_n}{(2k-1)/d_n}\gamma_n^{(2k-1)/d_n}\nonumber\\
&&\times\zeta_{n-1, \d', \gammab'}\left(N_1,\dots,N_{n-2}, N_{n-1}+N_n +(2k-1)/d_n\right).
\end{eqnarray*}
Lastly noting the fact that ${-N_n \choose (2k-1)/d_n}=0$ if $k>(1-d_n N_n)/2$,
we conclude the proof of Theorem \ref{main-main}.\qed

Next we proceed to the proofs of corollaries.

{\bf Proof of Corollary \ref{main1}:}
First prove point 1 of Corollary \ref{main1} by induction on $n$.
When $n=1$ it is clear, because
\begin{align}\label{cor1pt1n=1}
\zeta_{1,d_1,\gamma_1}(-N_1)=\gamma_1^{N_1}\zeta(-d_1 N_1)
=\gamma_1^{N_1}(-1)^{d_1 N_1}\frac{B_{d_1N_1+1}}{d_1N_1+1}.
\end{align}
The general case then follows from
the identity (\ref{mainmainFormula}) of Theorem \ref{main-main}  since 
$$\left(-N_1,\dots,-N_{n-2}, -N_{n-1}-N_n +(2k-1)/d_n\right)\in (-\N_0)^{n-1}$$ 
for $\Nb =(N_1,\dots,N_n) \in \N_0^n$ and $k\in\{1,\dots, \left[\frac{1+d_n N_n}{2}\right]\}$ such that $d_n|2 k-1$.\par

Assume now that $d_2,\dots,d_n$ are even integers satisfying the assumption (\ref{djsassumption}).
Let $\Nb=(N_1,\dots,N_n) \in \N_0^n$. 
Since now there is no $k$ for which $d_n|2k-1$ holds,
the identity  (\ref{mainmainFormula}) of Theorem \ref{main-main} implies that 
$$\zeta_{n, \d, \gammab} (-\Nb)= -\frac{1}{2}\zeta_{n-1, \d', \gammab'}\left(-N_1,\dots,-N_{n-2}, -N_{n-1}-N_n\right).$$
By induction on $n$, noting \eqref{cor1pt1n=1}, we deduce that 
\begin{eqnarray*}
\zeta_{n, \d, \gammab} (-\Nb)= \left(-\frac{1}{2}\right)^{n-1}\zeta_{1, d_1, \gamma_1} (-|\Nb|)
=\left(-\frac{1}{2}\right)^{n-1} \gamma_1^{|\Nb|}(-1)^{d_1|\Nb|}\frac{B_{d_1|\Nb|+1}}{(d_1|\Nb|+1)}.
\end{eqnarray*}
This ends the proof of point 2 of Corollary  \ref{main1}.
Lastly, point 3 follows immediately from point 2 because $B_1=-1/2$ and
$B_{2m+1}=0$ for all $m\in\mathbb{N}$.
\qed

{\bf Proof of Corollary \ref{main2}:} 
First, when $\Nb=\zerob_n$, from \eqref{mainmainFormula} we have
$$
\zeta_{n,\d,\gammab}(\zerob_n)=-\frac{1}{2}\zeta_{n-1,\d',\gammab'}(\zerob_{n-1}),
$$
from which point 1 immediately follows.

Next consider the case $\Nb=(\zerob_{n-1},-1)$.
The condition $d_n|2k-1$ implies whether $d_n=2k-1$ or $(2k-1)/d_n\geq 2$.
But in the latter case $\binom{-N_n}{(2k-1)/d_n}=\binom{1}{(2k-1)/d_n}=0$, so the
only possibility of $k$ on the right-hand side of \eqref{mainmainFormula} is
$k=(d_n+1)/2$.    When $d_n$ is odd this is indeed possible, and 
\eqref{mainmainFormula} gives
\begin{align*}
\zeta_{n,d,\gammab}(\zerob_{n-1},-1)&=-\frac{1}{2}\zeta_{n-1,\d',\gammab'}
(\zerob_{n-2},-1)-
\frac{B_{d_n+1}}{d_n+1}\gamma_n\zeta_{n-1,\d',\gammab'}(\zerob_{n-1})\\
&=-\frac{1}{2}\zeta_{n-1,\d',\gammab'}
(\zerob_{n-2},-1)-
\frac{B_{d_n+1}}{d_n+1}\gamma_n\left(-\frac{1}{2}\right)^{n-1}
\end{align*}
by using the result of point 1.   This formula is also valid for even $d_n$, because 
in this case $B_{d_n+1}=0$.    Using the above formula repeatedly, we obtain
$$
\zeta_{n,d,\gammab}(\zerob_{n-1},-1)=\left(-\frac{1}{2}\right)^{n-1}
\zeta_{1,d_1,\gamma_1}(-1)-\left(-\frac{1}{2}\right)^{n-1}
\sum_{j=2}^n \frac{B_{d_j+1}}{d_j+1}\gamma_j.
$$
Since $\zeta_{1,d_1,\gamma_1}(-1)=\gamma_1(-1)^{d_1}B_{d_1+1}/(d_1+1)$
(see \eqref{cor1pt1n=1}), we obtain
the assertion of point 2.

The case $\Nb=(\zerob_{n-1},-2)$ is similar.    In this case $N_n=-2$, so
$k\leq [(1-d_nN_n)/2]=[d_n+1/2]=d_n$.   Therefore $(2k-1)/d_n\geq 2$ is 
impossible, so the only possible $k$ is again $k=(d_n+1)/2$.    We obtain
\begin{align*}
\zeta_{n,d,\gammab}(\zerob_{n-1},-2)&=-\frac{1}{2}\zeta_{n-1,\d',\gammab'}
(\zerob_{n-2},-2)-
2\frac{B_{d_n+1}}{d_n+1}\gamma_n\zeta_{n-1,\d',\gammab'}(\zerob_{n-2},-1).
\end{align*}
We rewrite the term $\zeta_{n-1,\d',\gammab'}(\zerob_{n-2},-1)$ on the right-hand side
by using point 2, and then use the resulting formula repeatedly.
Lastly we use $\zeta_{1,d_1,\gamma_1}(-2)=\gamma_1^2B_{2d_1+1}/(2d_1+1)$ to arrive 
at the assertion of point 3.\qed

{\bf Remark:}
If we consider the case $\Nb=(\zerob_{n-1},-m)$, $m\geq 3$, larger values of $k$ 
appear on the right-hand side of \eqref{mainmainFormula}, so the explicit formula for
$\zeta_{n,d,\gammab}(\zerob_{n-1},-m)$ is (possible to obtain but) more complicated.
Therefore we only state the formula for $m\leq 2$ in Corollary \ref{main2}.
\medskip

{\bf Proof of Corollary \ref{intriguing}:}
Let $d_1\in \N\setminus\{1\}$ and $d_2\in \N$ be such that $\frac{1}{d_1}+\frac{1}{2d_2}\not \in \N$.\\
It follows that $(d_1, 2 d_2)$ satisfies the assumption (\ref{djsassumption}). Theorem \ref{main-main} implies then that 
for $N\in \N_0$, 
$$\zeta_{2,(d_1,2 d_2),(1,1)}(1+N, -N)=-\frac{1}{2}\zeta_{1,d_1,1}(1)=-\frac{1}{2}\zeta (d_1),$$
because $2d_2|(2k-1)$ is impossible.
This ends the proof of Corollary \ref{intriguing}.\qed

\section{Proof of Theorem \ref{irrationalityTh}}\label{sec-proofirrationalityTh}

In this section we assume (\ref{ira-assumption}).
Fix $\Nb \in \N_0^n$ and  $K\in \N$ such that $$2K >d_n(\frac{1}{d_1}+\dots+\frac{1}{d_n}) +d_n (N_1+\dots+ N_n ).$$
It is easy to see that $-\Nb \in \mathcal D_{n,\d}(-2K)$. Furthermore, Proposition \ref{mainproposition1} implies that there exists a function $\s\mapsto \mathcal R_{n,\d, K}(\s)$ holomorphic in 
$\mathcal D_{n,\d}(-2K)$ and satisfying $\mathcal R_{n,\d, K}(-\Nb) =0$ such that 
for all $\s \in \mathcal D_{n,\d}(0)$, we have 
\begin{eqnarray}\label{mainformula1bis}
\zeta_{n, \d, \gammab} (\s)&=&\frac{\Gamma(s_n-1/d_n) \Gamma(1/d_n) \gamma_n^{-1/d_n}}{d_n \Gamma (s_n)}\zeta_{n-1, \d', \gammab'}\left(s_1,\dots,s_{n-2}, s_{n-1}+s_n -1/d_n\right)\nonumber \\
&& +G_{n,\d,\gammab, K}(\s)+\mathcal R_{n,\d, K}(\s),
\end{eqnarray}
where
\begin{align}\label{G_def}
G_{n,\d,\gammab, K}(\s) &:=   -\frac{1}{2}\zeta_{n-1, \d', \gammab'}\left(s_1,\dots,s_{n-2}, s_{n-1}+s_n\right) \\
& -\sum_{k=1\atop d_n\mid 2k-1}^K \frac{B_{2k}}{2k}{-s_n \choose (2k-1)/d_n}\gamma_n^{(2k-1)/d_n}\notag\\
&\times\zeta_{n-1, \d', \gammab'}\left(s_1,\dots,s_{n-2}, s_{n-1}+s_n +(2k-1)/d_n\right).
\notag
\end{align}

Proposition \ref{basicproperties1} implies that $\s\mapsto G_{n,\d,\gammab, K}(\s)$ is meromorphic in $\C^n$. 
Moreover, if $\s=-\Nb$ is a singular point of $G_{n,\d,\gammab, K}(\s)$, point 2 of 
Proposition \ref{basicproperties1} implies then that there exists $j\in \{1,\dots,n-1\}$, $M, k_j\in \N_0$ and  $\eps_{j+1},\dots, \eps_{n-1} \in \{0,1\}$ such that 
$$-N_j-\dots- N_n +M=\frac{1}{d_j}+\frac{\eps_{j+1}}{d_{j+1}}+\dots +\frac{\eps_{n-1}}{d_{n-1}} -k_j.$$
It follows that $\frac{1}{d_j}+\frac{\eps_{j+1}}{d_{j+1}}+\dots +\frac{\eps_{n-1}}{d_{n-1}} \in \N$ which contradicts the assumption (\ref{ira-assumption}). 
As a conclusion, we find that
$\s=-\Nb$ is a regular point of $G_{n,\d,\gammab, K}$ and 
\begin{eqnarray}\label{ratpart}
G_{n,\d,\gammab, K}(-\Nb)&:=&-\frac{1}{2}\zeta_{n-1, \d', \gammab'}(-N_1,\dots,-N_{n-2}, -N_{n-1}-N_n)  \\
&& -\sum_{k=1\atop d_n \mid 2k-1}^{\left[(1+d_nN_n)/2\right]} \frac{B_{2k}}{2k} {N_n\choose (2k-1)/d_n} \gamma_n^{(2k-1)/d_n}\nonumber\\
&& \times  \zeta_{n-1, \d', \gammab'}\left(-N_1,\dots,-N_{n-2}, -N_{n-1}-N_n+(2k-1)/d_n\right).\nonumber
\end{eqnarray}

Let $k\in \{1,\dots,n-1\}$, and define
$$H_{n,\d,k}(\s):=\frac{\prod_{j=n-k+1}^n  \Gamma\left(s_j+\dots+s_n-\frac{1}{d_j}-\dots-\frac{1}{d_{n}}\right)
\Gamma\left(\frac{1}{d_j}\right) }{\prod_{j=n-k+1}^{n-1}  \Gamma\left(s_j+\dots+s_n-\frac{1}{d_{j+1}}-\dots-\frac{1}{d_{n}}\right) {\prod_{j=n-k+1}^{n}d_j \gamma_j^{1/d_j}}}.$$
We now prove the following
\begin{lemma}\label{interm_lemma_theorem2}
For $k\in \{1,\dots,n-1\}$, 
there exists a function $\s\mapsto \mathcal R_{n,\d, K,k}(\s)$, meromorphic in 
$\mathcal D_{n,\d}(-2K)$, regular at $\s=-\Nb$ and satisfying $\mathcal R_{n,\d, K,k}(-\Nb) =0$,
for which the formula
\begin{align}\label{big-induction}
\lefteqn{\zeta_{n, \d, \gammab} (\s)=\frac{1}{\Gamma(s_n)}H_{n,\d,k}(\s)}\nonumber\\
&  \;\;\times \zeta_{n-k, (d_1,\dots,d_{n-k}), (\gamma_1,\dots,\gamma_{n-k})}\left(s_1,\dots,s_{n-k-1}, s_{n-k}+\dots+s_n -\frac{1}{d_{n-k+1}}-\dots-\frac{1}{d_{n}}\right) \nonumber\\
& +G_{n,\d,\gammab, K}(\s)+ \mathcal R_{n,\d, K,k}(\s)
\end{align}
holds for all $\s\in  \mathcal D_{n,\d}(-2K)$ as an identity of meromorphic functions, 
where $G_{n,\d,\gammab, K}$ is the function defined by (\ref{G_def}).
\end{lemma}

{\bf Proof:}
The proof is by induction on $k$.
For $k=1$, (\ref{big-induction}) holds from (\ref{mainformula1bis}) and (\ref{ratpart}).

Let $k\in \{1,\dots, n-2\}$. Assume that (\ref{big-induction}) holds for $k$. We will prove that it also holds for $k+1$.

Let $K'\in \N$ be such that $K'>\max \left(K, K\frac{d_{n-k}}{d_n}\right)$.
Since 
$$\left(s_1,\dots,s_{n-k-1}, s_{n-k}+\dots+s_n -\frac{1}{d_{n-k+1}}-\dots-\frac{1}{d_{n}}\right) \in  \mathcal D_{n-k,(d_1,\dots,d_k)}(0)$$ 
for all $\s\in  \mathcal D_{n,\d}(0)$,
applying Proposition \ref{mainproposition1} (with $K'$ instead of $K$) to the term 
$$ \zeta_{n-k, (d_1,\dots,d_{n-k}), (\gamma_1,\dots,\gamma_{n-k})}\left(s_1,\dots,s_{n-k-1}, s_{n-k}+\dots+s_n -\frac{1}{d_{n-k+1}}-\dots-\frac{1}{d_{n}}\right)$$ 
on the right-hand side of (\ref{big-induction}), for $\s \in \mathcal D_{n,\d}(0)$, we obtain
\begin{align*}
\lefteqn{\zeta_{n, \d, \gammab} (\s)=\frac{1}{\Gamma(s_n)}H_{n,\d,k+1}(\s)}\\
& \;\;\times \zeta_{n-k-1, (d_1,\dots,d_{n-k-1}), (\gamma_1,\dots,\gamma_{n-k-1})}\left(s_1,\dots,s_{n-k-2}, s_{n-k-1}+\dots+s_n -\frac{1}{d_{n-k}}-\dots-\frac{1}{d_{n}}\right) \\
& +G_{n,\d,\gammab, K}(\s)+ \mathcal R_{n,\d, K,k}(\s) +\frac{1}{\Gamma(s_n)} 
H_{n,\d,k}(\s) V_{n,\d,k}(\s),
\end{align*}
where
\begin{align*}
\lefteqn{V_{n,\d,k}(\s):=
-\frac{1}{2}  \zeta_{n-k-1, (d_1,\dots,d_{n-k-1}), (\gamma_1,\dots,\gamma_{n-k-1})}\biggl(s_1,\dots,s_{n-k-2}, s_{n-k-1}+\dots+s_n\biggr.}\\
&\qquad\qquad\qquad\biggl. -\frac{1}{d_{n-k+1}}-\dots-\frac{1}{d_{n}}\biggr) \\
&  -\sum_{k=1\atop d_{n-k}\mid 2k-1}^{K'} \frac{B_{2k}}{2k}{-s_{n-k}-\dots-s_n+1/d_{n-k+1}+\dots+1/d_n \choose (2k-1)/d_{n-k}}\gamma_{n-k}^{(2k-1)/d_{n-k}} \\
& \qquad\times \zeta_{n-k-1, (d_1,\dots,d_{n-k-1}), (\gamma_1,\dots,\gamma_{n-k-1})}\biggl(s_1,\dots,s_{n-k-2}, s_{n-k-1}+\dots+s_n \biggr.\\
&\qquad\biggl.-\frac{1}{d_{n-k+1}}-\dots-\frac{1}{d_{n}}+\frac{2k-1}{d_{n-k}}\biggr)\\
&+\mathcal R_{n-k, (d_1,\dots,d_{n-k}),K'}\left(s_1,\dots,s_{n-k-1}, s_{n-k}+\dots+s_n -\frac{1}{d_{n-k+1}}-\dots-\frac{1}{d_{n}}\right).
\end{align*}
Since $K'>\max \left(K, K\frac{d_{n-k}}{d_n}\right)$, it is easy to see that the function 
$$\s\mapsto \mathcal R_{n-k, (d_1,\dots,d_{n-k}),K'}\left(s_1,\dots,s_{n-k-1}, s_{n-k}+\dots+s_n -\frac{1}{d_{n-k+1}}-\dots-\frac{1}{d_{n}}\right)$$ 
is holomorphic in $\mathcal D_{n,\d}(-2K)$, it follows then from  Proposition \ref{basicproperties1} that $\s\mapsto V_{n,\d,k}(\s)$ is a meromorphic function in $\mathcal D_{n,\d}(-2K)$. 
Moreover, if $\s=-\Nb$ is a singular point of $V_{n,\d,k}(\s)$, point 3 of Proposition \ref{basicproperties1} implies then that there exists $j\in \{1,\dots,n-k-1\}$, $M, k_j\in \N_0$ and  $\eps_{j+1},\dots, \eps_{n-k-1} \in \{0,1\}$ such that 
$$-N_j-\dots- N_n -\frac{1}{d_{n-k+1}}-\dots -\frac{1}{d_{n}}+M=\frac{1}{d_j}+\frac{\eps_{j+1}}{d_{j+1}}+\dots+\frac{\eps_{n-k-1}}{d_{n-k-1}} -k_j.$$
It follows that $\frac{1}{d_j}+\frac{\eps_{j+1}}{d_{j+1}}+\dots+\frac{\eps_{n-k-1}}{d_{n-k-1}}+\frac{1}{d_{n-k+1}}+\dots +\frac{1}{d_{n}} \in \N$ which contradicts the assumption (\ref{ira-assumption}). 
As a conclusion, we prove that 
$\s\mapsto V_{n,\d,k}(\s)$ is regular at $\s=-\Nb$.

Also, assumption (\ref{ira-assumption}) implies that $\s\mapsto H_{n,\d,k}(\s)$ is a meromorphic function in  
$\C^n$  which is regular at 
$\s=-\Nb$.   By using in addition the fact that $\frac{1}{\Gamma (s_n)}|_{s_n=-N_n}=0$, we deduce that 
$$s\mapsto \mathcal R_{n,\d, K,k+1}(\s):=\mathcal R_{n,\d, K,k}(\s) +\frac{1}{\Gamma(s_n)} H_{n,\d,k}(\s) V_{n,\d,k}(\s)$$
 is is a meromorphic function in  $\mathcal D_{n,\d}(-2K)$  which is regular at $\s=-\Nb$ 
 and satisfies $R_{n,\d, K,k+1}(-\Nb)=0$. 
 This ends the induction argument and therefore ends the proof of (\ref{big-induction}).\qed
 
We are now ready to finish the proof of Theorem \ref{irrationalityTh}. It follows from  (\ref{big-induction}) with $k=n-1$ that 
\begin{eqnarray}\label{big-inductionbis}
\zeta_{n, \d, \gammab} (\s)&=&\frac{\prod_{j=2}^n\Gamma\left(\frac{1}{d_j}\right) }{\prod_{j=2}^n d_j\gamma_j^{1/d_j}}~
A(\s)~B(\s)~ \zeta_{1,d_1,\gamma_1}(s_1+\cdots+s_n-b) \nonumber\\
&& +G_{n,\d,\gammab, K}(\s)+ \mathcal R_{n,\d, K,n-1}(\s),
\end{eqnarray}
where  
\be
\item  $b:=\frac{1}{d_2}+\dots+\frac{1}{d_{n}}\in \N$;
\item $G_{n,\d,\gammab, K}$ is the function defined by (\ref{G_def}) and $\s\mapsto \mathcal R_{n,\d, K,n-1}(\s)$ is a meromorphic function in 
$\mathcal D_{n,\d}(-2K)$, regular in $\s=-\Nb$  and satisfying $\mathcal R_{n,\d, K,n-1}(-\Nb) =0$;
\item 
$\ds A(\s):=\frac{\prod_{j=3}^n  \Gamma\left(s_j+\dots+s_n-\frac{1}{d_j}-\dots-\frac{1}{d_{n}}\right)}{\prod_{j=2}^{n-1}  \Gamma\left(s_j+\dots+s_n-\frac{1}{d_{j+1}}-\dots-\frac{1}{d_{n}}\right) },\; B(\s):=\frac{\Gamma\left(s_2+\dots+s_n-b\right)}{\Gamma (s_n)}.$
\ee
Assumption (\ref{ira-assumption}) implies that $\s\mapsto A(\s)$ is regular in $\s=-\Nb$ and that
\begin{equation}\label{A-N}
A(-\Nb)=\frac{\prod_{j=3}^n  \Gamma\left(-N_j-\dots-N_n-\frac{1}{d_j}-\dots-\frac{1}{d_{n}}\right)}{\prod_{j=2}^{n-1}  \Gamma\left(-N_j-\dots-N_n-\frac{1}{d_{j+1}}-\dots-\frac{1}{d_{n}}\right) }.
\end{equation}
Moreover, by using the identity 
$\Gamma (z+1)=z \Gamma (z)$, it is easy to see that for $M \in \N_0$ and $x \in \C\setminus \Z$, we have 
\begin{equation}\label{gammadecalage}
\Gamma (-M-x)=\frac{\Gamma (b-x)}{(x)_{M,b}}\quad {\mbox{ where }}\quad (x)_{M,b}=\prod_{k=-M}^{b-1} (k-x).
\end{equation}
Combining (\ref{A-N}) and (\ref{gammadecalage}) we have
\begin{eqnarray}\label{A-N-final}
A(-\Nb)&=&\frac{\prod_{j=3}^n  \Gamma\left(b-\frac{1}{d_j}-\dots-\frac{1}{d_{n}}\right)}{\prod_{j=3}^n
\left(\frac{1}{d_j}+\dots+\frac{1}{d_{n}}\right)_{N_j+\dots+N_n, b}}
\frac{\prod_{j=2}^{n-1} \left(\frac{1}{d_{j+1}}+\dots+\frac{1}{d_{n}}\right)_{N_j+\dots+N_n, b}}{\prod_{j=2}^{n-1}
  \Gamma\left(b-\frac{1}{d_{j+1}}-\dots-\frac{1}{d_{n}}\right)}\nonumber\\
&=&\frac{\prod_{j=2}^{n-1} \left(\frac{1}{d_{j+1}}+\dots+\frac{1}{d_{n}}\right)_{N_j+\dots+N_n, b}}{\prod_{j=3}^n
\left(\frac{1}{d_j}+\dots+\frac{1}{d_{n}}\right)_{N_j+\dots+N_n, b}}
=\prod_{j=3}^n\frac{\left(\frac{1}{d_j}+\dots+\frac{1}{d_{n}}\right)_{N_{j-1}+\dots+N_n, b}}{\left(\frac{1}{d_j}+\dots+\frac{1}{d_{n}}\right)_{N_j+\dots+N_n, b}}\nonumber\\
&=&\prod_{j=3}^n \left(\prod_{u=-N_{j-1}-\dots-N_n}^{-N_{j}-\dots-N_n-1}\left(u-\sum_{k=j}^n \frac{1}{d_k}\right)\right).
\end{eqnarray}

On the other hand, $\s=-\Nb$ is a singular point of $\s\mapsto B(\s)$ and it is a point of indeterminacy. 
Fix $\thetab=(\theta_1,\dots,\theta_n) \in \C^n$ such that $\theta_2+\dots+\theta_n\neq 0$ and $\theta_n\neq 0$. Set 
$\delta = \min \left\{|\theta_2+\dots+\theta_n|^{-1}, |\theta_n|^{-1}\right\} >0$. Then, for all $t\in \C\setminus\{0\}$ such that $|t|<\delta$, we have 
$$B(-\Nb+t\thetab)=\frac{\Gamma\left(-(N_2+\dots+N_n+b)+t(\theta_2+\dots+\theta_n)\right)}{\Gamma (-N_n+t\theta_n)}.$$
By using the classical fact that for $k\in \N_0$, $\Gamma(z)$ has a simple pole in $z=-k$ of residue $\frac{(-1)^k}{k!}$, we deduce that the  directional limit $B^{\thetab}(-\Nb):=\lim_{t\rightarrow 0, ~t\in \C\setminus\{0\}} B(-\Nb+t\thetab)$ exists and is given by 
\begin{equation}\label{B-Nfinal}
B^{\thetab}(-\Nb) = \frac{(-1)^{N_2+\dots+N_{n-1}+b} N_n!}{(N_2+\dots+N_n+b)!}\cdot \left(\frac{\theta_n}{\theta_2+\dots+\theta_n}\right).
\end{equation}
Combining (\ref{big-inductionbis}), (\ref{A-N}), (\ref{B-Nfinal}) and (\ref{ratpart}) 
we find that 
the directional limit $\zeta_{n, \d, \gammab}^{\thetab}(-\Nb):=
\lim_{t\rightarrow 0, ~t\in \C\setminus\{0\}}\zeta_{n, \d, \gammab} (\Nb+t\theta)$ exists and is given by 
$$
\zeta_{n, \d, \gammab}^{\thetab} (-\Nb)=\frac{\prod_{j=2}^n\Gamma\left(\frac{1}{d_j}\right) }{\prod_{j=2}^n d_j\gamma_j^{1/d_j}}~
A(-\Nb)~B^{\thetab}(-\Nb)~\zeta_{1,d_1,\gamma_1}(-(|\Nb| +b)) 
 +G_{n,\d,\gammab, K}(-\Nb).
$$
We conclude by using in addition the expressions given by (\ref{cor1pt1n=1}), (\ref{ratpart}), (\ref{A-N-final}) and (\ref{B-Nfinal}). This ends the proof of  Theorem \ref{irrationalityTh}.\qed

\section{Values of Mahler's series at non-positive integers}\label{sec-mahlerseries}

Now we proceed to the discussion of more general zeta-functions
$\zeta_n(\s;\P)$.
In this section, as a preparation, we study the values of multiple series of Mahler type at
non-positive integer points. We will use notations introduced in the beginning of Section 3. 

Let $P\in \R[X_1,\dots, X_n]$ be an elliptic polynomial of degree $d\geq 1$ and let  $Q\in \R[X_1,\dots, X_n]$ be a polynomial of degree $q\geq 0$.
Define $\mathcal D:=\{\Re(s) >\frac{n+q}{d}\}$,
$$Y(P,Q;s):= 
\int_{[1,\infty)^n} Q (\x) P^{-s}(\x) ~d\x \quad {\mbox { and }}\quad  Z(P,Q;s):= 
\sum_{\m \in \N^n} \frac{Q(\m)}{ P(\m)^{s}}. $$ 

\begin{lemma}\label{lem_mahler}
{\rm (K.~Mahler \cite{mahler})}
Both $Y(P,Q;s)$ and 
$Z(P, Q; s)$ converges absolutely in 
$\mathcal D$ and has meromorphic continuation to the whole complex plane $\C$ with at most simple poles located in the set 
$$\mathcal P (P,Q):=\left\{s=\frac{n+q-k}{d}\mid k \in \N_0\right\}\setminus (-\N_0).$$
\end{lemma}

For elliptic polynomials $P$, Pierrette Cassou-Nogu\`es obtained  in the eighties  (\cite{CN1}, \cite{CN2}, etc.)  several important results on the values of $ Z(P, Q; -N)$ at non-positive integers  $-N$.
(cf. \cite{MW} for another approach.)

In the proof of Theorem \ref{mainsection3}, we will use the following result which gives new closed formulas  for the values of $ Z(P, Q; -N)$. Our proof of this result  is also different from the method of
Cassou-Nogu\`es.  

\begin{theorem}\label{mainbissection3}
Let $P\in \R[X_1,\dots, X_n]$ be an elliptic and homogeneous polynomial of degree $d\geq 1$ and let  $Q\in \R[X_1,\dots, X_n]$ be a homogeneous polynomial of degree $q\geq 0$.
Then, for any $N\in \N_0$, $s=-N$ is not a pole of $Z(P,Q;s)$ and 
 \begin{eqnarray*}
 Z(P,Q; -N) &=&
 \sum_{\betab \in \N_0^n \atop |\betab|\leq q} \sum_{\alphab \in I_N(\betab)} \sum_{\u \in V(\alphab)}
 \frac{(-1)^{|\alphab|-N} (|\alphab|-1-N)! N!}{d ~\alphab ! ~\betab !} \\
 && \quad \times \left(\sum_{i=1}^n K_i(P;Q; N; \alphab, \u, \betab)\right) ~~
 \prod_{i=1}^n {\widetilde{B}}_{g_i(\u)+\beta_i},
 \end{eqnarray*}
where ${\widetilde{B}}_k :=B_k$ is as in the statement of Theorem \ref{mainsection3}.
\end{theorem}

{\bf Proof of Theorem \ref{mainbissection3}:}\\
Let $P\in \R[X_1,\dots, X_n]$ be an elliptic and homogeneous polynomial of degree $d\geq 1$ and let  
$Q\in \R[X_1,\dots, X_n]$ be a homogeneous polynomial of degree $q\geq 0$.
Define for any $\a=(a_1,\dots,a_n) \in [0,\infty)^n$
$$P_\a =P(\X+\a) =P(X_1+a_1,\dots, X_n+a_n) {\mbox { and }} Q_\a =Q(\X+\a) =Q(X_1+a_1,\dots, X_n+a_n).$$
Define for any 
$s \in \mathcal D$ and any $\a=(a_1,\dots,a_n) \in [0,\infty)^n$
$$Y(P_\a,Q_\a;s)= 
\int_{[1,\infty)^n} Q_\a (\x) P_\a^{-s}(\x) ~d\x \quad {\mbox { and }}\quad  Z(P_\a,Q_\a;s)= 
\sum_{\m \in \N^n} \frac{Q_\a (\m)}{ P_\a(\m)^{s}}. $$ 
Our first useful ingredient is the following result:

\begin{prop}\label{ypaqas}
Let $\a \in [0,\infty)^n$.   The integral
$Y(P_\a,Q_\a;s)$ converges absolutely in  $\mathcal D:=\{ \Re(s) >\frac{n+q}{d}\}$ and has a meromorphic continuation to $\C$ with at most simple poles located in the set 
$\mathcal P (P,Q)$ defined in Lemma \ref{lem_mahler}. Moreover, 
for any $N\in \N_0$,  $Y(P_\a,Q_\a;s)$ is regular at $s=-N$ and its value is given by 
\begin{eqnarray*}
 Y(P_\a,Q_\a; -N) &=&
 \sum_{\betab \in \N_0^n \atop |\betab|\leq q} \sum_{\alphab\in I_N(\betab)} \sum_{\u \in V(\alphab)}
 \frac{(-1)^{|\alphab|-N} (|\alphab|-1-N)! N!}{d ~\alphab ! ~\betab !} \\
 && \quad \times \left(\sum_{i=1}^n K_i(P;Q; N; \alphab, \u, \betab)\right) ~~
 \prod_{i=1}^n (1+a_i)^{g_i(\u)+\beta_i}.
 \end{eqnarray*}
\end{prop}

{\bf Proof of Proposition \ref{ypaqas}:}\\
Fix $\a \in [0,\infty)^n$ and set $\b =(b_1,\dots, b_n)$ where $b_i =a_i +1$ 
($i=1,\dots,n$). 
First we remark that there exist $C_1, C_2$ and  $C_3 >0$ such that 
$$|Q(\x)|\leq C_1 |\x|^q {\mbox { and }} C_2 |\x|^d \leq P(\x) \leq C_3 |\x|^d \quad
\mbox{for all}\;\; \x \in [1,\infty)^n.$$
Then, since
$Q_\a (\x) P_\a^{-s}(\x)\ll |\x|^{q-d\Re s}$ and if $\Re s > (n+q)/d$ then
$q-d\Re s < -n$, we see that
$Y(P_\a,Q_\a;s)$ converges absolutely in  $\mathcal D$.

Set for all $i=1,\dots,n$,
$$V_i=\big\{\x\in (0,\infty)^n \mid x_j < x_i \;\mbox{for all}\; j\in \{1,\dots,n\}\setminus \{i\}\big\}.$$
Shifting $\x$ to $\x+1$ in the definition integral of $Y(P_\a,Q_\a;s)$, we find that
for any $s\in \mathcal D$:
\begin{equation}\label{decoforsing}
Y(P_\a,Q_\a;s)=\sum_{i=1}^n Y_i(\b; s) \quad {\mbox { where }} \quad  Y_i(\b; s):=\int_{V_i} Q_\b (\x) P_\b^{-s}(\x) ~d\x.
\end{equation}
We also use the notation, for all $i=1,\ldots,n$,
$$H_k\left({\hat{\y}}(i);\b\right)= \sum_{\gammab \in \N_0^n \atop |\gammab|=k} \frac{\b^{\gammab}}{\gammab!} 
\partial^{\gammab}P \left({\hat{\y}}(i))\right)\quad \;\mbox{for all}\; k=1,\dots,d$$
and
$$\H\left({\hat{\y}}(i);\b\right)=\left(H_1\left({\hat{\y}}(i);\b\right),\dots, H_d\left({\hat{\y}}(i);\b\right)\right).$$

We will first study $Y_n(\b; s)$.
The Taylor expansion implies that for any $s\in \mathcal D$:
\begin{equation}\label{decoforsing2}
Y_n(\b;s)=\sum_{\betab \in \N_0^n \atop |\betab|\leq q} \frac{\b^{\betab}}{\betab!} Y_n(\b; \betab; s) \quad {\mbox { where }} \quad  
Y_n(\b; \betab; s):=\int_{V_n} (\partial^{\betab}Q (\x)) P_\b^{-s}(\x) ~d\x,
\end{equation}
because $\partial^{\betab}Q (\x)=0$ if $|\betab|> q$.

Fix $\betab \in \N_0^n$ such that $|\betab|\leq q$.
Consider the blowing up 
$$\varphi_n : (0,1)^{n-1}\times (0,\infty) \rightarrow V_n,$$
defined by
$$\y=(y_1,\dots, y_n)\mapsto \varphi_n (\y)=(y_1y_n,\dots, y_{n-1}y_n, y_n).$$
Since this is a map onto $V_n$ and its Jacobian is $y_n^{n-1}$, 
we get that for any $s\in \mathcal D$:
$$Y_n(\b; \betab; s)=\int_{(0,1)^{n-1}\times (0,\infty)} (\partial^{\betab}Q (\varphi_n(\y))) P_\b^{-s}(\varphi_n(\y)) ~y_n^{n-1}~d\y.$$
(Here, it is to be noted that $\partial^{\betab}Q (\varphi_n(\y))$ does not mean the
derivative in $\y$; it means to substitute $\varphi_n(\y)$ in place of $\x$ into 
$\partial^{\betab}Q (\x)$.)
Since $\partial^{\betab}Q (\x)$ is also homogeneous, we obtain
\begin{eqnarray}\label{blowup}
Y_n(\b; \betab; s)=
 \int_{(0,1)^{n-1}\times (0,\infty)} (\partial^{\betab}Q ({\hat{\y}}(n))) P_\b^{-s}(\varphi_n(\y)) ~y_n^{q-|\betab|+n-1}~d\y.
\end{eqnarray}
Again applying the Taylor formula we see that for any $\y \in (0,1)^{n-1}\times (0,\infty)$
\begin{eqnarray}\label{pnphiybexp}
P_\b\left(\varphi_n(\y)\right)&=& P\left(\varphi_n(\y)+\b\right)=\sum_{\gammab \in \N_0^n \atop |\gammab|\leq d} \frac{\b^{\gammab}}{\gammab!}
\partial^{\gammab}P \left(\varphi_n(\y)\right)
=\sum_{\gammab \in \N_0^n \atop |\gammab|\leq d} \frac{\b^{\gammab}}{\gammab!} y_n^{d-|\gammab|}
\partial^{\gammab}P \left({\hat{\y}}(n)\right)\nonumber\\
&=&y_n^d P \left({\hat{\y}}(n)\right) + \sum_{k=1}^d y_n^{d-k} H_k\left({\hat{\y}}(n);\b\right). 
\end{eqnarray}
The ellipticity of $P$ imply that 
\begin{equation}\label{ellipticity}
P \left({\hat{\y}}(n))\right) = P(y_1,\dots, y_{n-1},1)>0 \quad \mbox{for all}\; \;(y_1,\dots, y_{n-1}) \in [0,1]^{n-1}.
\end{equation}
Therefore, we deduce from the compacity of $[0,1]^{n-1}$ that 
\begin{equation}\label{boundbien}
B=B(P;\b):=\sup_{(y_1,\dots,y_{n-1})\in [0,1]^{n-1}} \sum_{k=1}^d \left|\frac{H_k\left({\hat{\y}}(n);\b\right)}{P \left({\hat{\y}}(n)\right)}\right|<\infty.
\end{equation}

Set $A=A(P;\b)= 2 B(P; \b)+2 \geq 2$, and divide the integral on the right-hand side 
of \eqref{blowup} as
\begin{align}\label{divide-two-integrals}
Y_n(\b; \betab; s)&=
 \int_{(0,1)^{n-1}\times (A,\infty)} +\int_{(0,1)^{n-1}\times (0,A)}\notag\\
 &=Y_n^{A+}(\b; \betab; s)+Y_n^{A-}(\b; \betab; s),\quad\mbox{say.}
\end{align}

We first consider the integral $Y_n^{A+}(\b; \betab; s)$.
For all $\y=(y_1,\dots, y_n) \in [0,1]^{n-1}\times [A, \infty)$, it follows that
\begin{equation}\label{estimateforconv}
\left|\sum_{k=1}^d \frac{H_k\left({\hat{\y}}(n);\b\right)}{y_n^k~P \left({\hat{\y}}(n)\right)}\right| \leq \frac{B}{y_n} \leq \frac{B}{A} < \frac{1}{2}.
\end{equation}
Now recall the elementary identity
\begin{equation}\label{elementary}
(1+X_1+\cdots+X_d)^{-s}=\sum_{\alphab \in \N_0^d}{-s\choose |\alphab|} \frac{|\alphab|!}{\alphab !}X_1^{\alpha_1}\cdots X_d^{\alpha_d}
\quad (|X_1+\cdots+X_d|<1).
\end{equation}
Applying this to
(\ref{pnphiybexp}) and using the upper bound (\ref{estimateforconv}) we find that for any compact subset $K$ of $\C$, 
we have 
\begin{align}\label{Pb}
P_\b\left(\varphi_n(\y)\right)^{-s}= \sum_{\alphab \in \N_0^d} {-s\choose |\alphab|} \frac{|\alphab|!}{\alphab !}\H\left({\hat{\y}}(n);\b\right)^{\alphab} 
P\left({\hat{\y}}(n)\right)^{-s-|\alphab|} y_n^{-d s-\sum_{k=1}^d k \alpha_k},
\end{align}
whose convergence is uniform in $s\in K$ and in $\y \in [0,1]^{n-1}\times [A, \infty)$.
The uniformity of the  convergence implies that for any 
$\s\in \mathcal D =\{\Re (s)> \frac{q+n}{d}\}$,
we can substitute \eqref{Pb} into the definition of $Y_n^{A+}(\b; \betab; s)$ and
carry out the termwise integration.    We then obtain
\begin{eqnarray}\label{blowup+}
Y_n^{A+}(\b; \betab; s)
= \sum_{\alphab \in \N_0^d} {-s\choose |\alphab|} \frac{|\alphab|!}{\alphab !}
\frac{A^{-d s-\sum_{k=1}^d k \alpha_k+q+n-|\betab|}}{d s+\sum_{k=1}^d k \alpha_k-q-n+|\betab|}
R(\alphab; s),
\end{eqnarray}
where
$$R(\alphab; s):=\int_{(0,1)^{n-1}} P\left({\hat{\y}}(n))\right)^{-s-|\alphab|} \partial^{\betab}Q \left({\hat{\y}}(n))\right)
\H\left({\hat{\y}}(n);\b\right)^{\alphab} ~dy_1\dots dy_{n-1}.$$
The function $s\mapsto R(\alphab; s)$ 
is clearly holomorphic in the whole complex plane $\C$. Moreover, the bound (\ref{boundbien}) implies $|H_k\left({\hat{\y}}(n);\b\right)|\leq B|P({\hat{\y}}(n))|$, so
for any compact subset $K$ of $\C$, 
we have uniformly in $\alphab \in \N_0^d$ and $s\in K$,
$$R(\alphab; s) \ll_{K,P,Q} B^{|\alphab|}< \left(\frac{A}{2}\right)^{|\alphab|}$$
(similar to \eqref{estimateforconv}).
Therefore, the last member of (\ref{blowup+}) defines a meromorphic function in $\C$ with at most simple poles located in the set 
$\mathcal P_0 (P,Q):=\{\frac{n+q-k}{d}\mid k\in \N_0\}$.
That is, $Y_n^{A+}(\b; \betab; s)$ has the meromorphic continuation (given explicitly by the last member of (\ref{blowup+})) to the whole complex plane $\C$ with at most simple poles located in the set 
$\mathcal P_0 (P,Q)$.

Now let $N\in \N_0$.   We consider the situation at $s=-N$.
The denominators of the terms in the sum (\ref{blowup+}) correspond to 
$\alphab \in I_N(\betab)$ are $0$ at
$s=-N$, so this point is a possible pole. 
However, if $\alphab \in  I_N(\betab)$, then $|\alphab|>N$ 
because
$$
dN<dN+n+(q-|\betab|)=\sum_{k=1}^d k\alpha_k\leq d\sum_{k=1}^d\alpha_k=d|\alphab|.
$$
Hence ${-s\choose |\alphab|}\big|_{s=-N}=0$. 
Therefore $s=-N$ is not a pole of $Y_n^{A+}(\b; \betab; s)$, and hence we conclude that
possible poles of $Y_n^{A+}(\b; \betab; s)$ are located only on $\mathcal P(P,Q)$. 

Using the fact that for any $\alphab \in \N_0^d$ such that $|\alphab|>N$ we have
$${-s\choose |\alphab|}\sim \frac{(-1)^{|\alphab|-N}N!(\alphab|-1-N)!}{|\alphab|!} ~(s+N) \quad {\mbox {as }} s\rightarrow -N,$$
we deduce then from (\ref{blowup+}) that 
\begin{eqnarray}\label{yna+bbeta-N}
Y_n^{A+}(\b; \betab; -N)&=&-\sum_{\alphab \in \N_0^d \atop |\alphab|\leq N} {N\choose |\alphab|} \frac{|\alphab|!}{\alphab !}
\frac{A^{d N-\sum_{k=1}^d k \alpha_k+q+n-|\betab|}}{d N+q+n-\sum_{k=1}^d k \alpha_k -|\betab|}
R(\alphab,-N)\nonumber\\
&& +\sum_{\alphab \in I_N(\betab)}  \frac{(-1)^{|\alphab|-N} N! (|\alphab|-1-N)!}{d ~\alphab !}R(\alphab,-N).
\end{eqnarray}

Next consider $Y_n^{A-}(\b; \betab; s)$.
Since $P$ is elliptic and $q\geq|\betab|$, we see that
$Y_n^{A-}(\b; \betab; s)$
is holomorphic in the whole complex plane $\C$, 
and for any $N\in \N_0$,
\begin{equation}\label{toto1}
Y_n^{A-}(\b; \betab; -N)
= \int_{(0,1)^{n-1}\times (0, A)} \partial^{\betab}Q \left({\hat{\y}}(n)\right) P_\b^{N}(\varphi_n(\y)) ~y_n^{q-|\betab|+n-1}~d\y.
\end{equation}
Furthermore, it follows from (\ref{pnphiybexp}) and \eqref{elementary} that 
\begin{eqnarray}\label{toto2}
P_\b^{N}(\varphi_n(\y)) &=& \left(y_n^d P \left({\hat{\y}}(n)\right) + \sum_{k=1}^d y_n^{d-k} H_k\left({\hat{\y}}(n);\b\right)\right)^N\nonumber\\
&=& \sum_{\alphab \in \N_0^d \atop |\alphab|\leq N} {N\choose |\alphab|}\frac{|\alphab|!}{\alphab !} y_n^{d N-\sum_{k=1}^d k \alpha_k}
P^{N-|\alphab|} \left({\hat{\y}}(n)\right) \H\left({\hat{\y}}(n);\b\right)^{\alphab}.
\end{eqnarray}
Combining (\ref{toto1}) and (\ref{toto2}) we obtain for any $N\in \N_0$,
\begin{eqnarray}\label{yna-bbeta-N}
\lefteqn{Y_n^{A-}(\b; \betab; -N)}\nonumber\\
&=&\sum_{\alphab \in \N_0^d \atop |\alphab|\leq N} {N\choose |\alphab|} \frac{|\alphab|!}{\alphab !}
\frac{A^{d N-\sum_{k=1}^d k \alpha_k+q+n-|\betab|}}{d N+q+n-\sum_{k=1}^d k \alpha_k -|\betab|}
R(\alphab,-N).
\end{eqnarray}
As a conclusion, (\ref{decoforsing2}), (\ref{yna+bbeta-N}) and (\ref{yna-bbeta-N}) imply that 
$$
Y_n(\b;s)=\sum_{\betab \in \N_0^n \atop |\betab|\leq q} \frac{\b^{\betab}}{\betab!} Y_n(\b; \betab; s) 
= \sum_{\betab \in \N_0^n \atop |\betab|\leq q} \frac{\b^{\betab}}{\betab!} \left(Y_n^{A-}(\b; \betab; s)+Y_n^{A-}(\b; \betab; s)\right)
$$
has a meromorphic  continuation to $\C$ with at most simple poles located on the set 
$\mathcal P (P,Q)$ and that for any $N\in \N_0$,
\begin{eqnarray}\label{yn+bbeta-N}
Y_n(\b; \betab; -N)&=&\sum_{\alphab \in I_N(\betab)}  \frac{(-1)^{|\alphab|-N} N! (|\alphab|-1-N)!}{d ~\alphab !}R(\alphab,-N),
\end{eqnarray}
because the first term on the right-hand side of \eqref{yna+bbeta-N} is cancelled with
the right-hand side of \eqref{yna-bbeta-N}.

By a simple permutation of the variables, we deduce from the previous argument that for any $i=1,\dots,n$, 
$Y_i(\b;s)$ has a meromorphic  continuation to $\C$ with at most simple poles located in the set 
$\mathcal P(P,Q)$ and that  for any $N\in \N_0$
\begin{eqnarray}\label{yi+bbeta-N}
Y_i(\b; \betab; -N)&=&\sum_{\alphab \in I_N(\betab)}  \frac{(-1)^{|\alphab|-N} N! (|\alphab|-1-N)!}{d ~\alphab !}\\
&& \times \int_{(0,1)^{n-1}} P\left({\hat{\y}}(i)\right)^{N-|\alphab|} (\partial^{\betab}Q \left({\hat{\y}}(i)\right))
\H\left({\hat{\y}}(i);\b\right)^{\alphab} ~\prod_{k=1\atop k\neq i}^n dy_k.\nonumber
\end{eqnarray}
To end the proof of Proposition \ref{ypaqas}, it suffices to remark that for any $i=1,\dots,n$:
\begin{eqnarray*}
\H\left({\hat{\y}}(i);\b\right)^{\alphab}&=&\prod_{k=1}^d \left(\sum_{\gammab \in \Delta_k^n} \frac{\b^{\gammab}}{\gammab!} 
\partial^{\gammab}P \left({\hat{\y}}(i))\right)\right)^{\alpha_k}\\
&=&\prod_{k=1}^d \left[\sum_{ \u_k =\left(u_{k,\gammab}\right)_{\gammab \in \Delta_k^n}\in \N_0^{{n+k-1 \choose n-1}},~ |\u_k|=\alpha_k}\frac{\alpha_k!}{\u_k !}
\prod_{\gammab \in \Delta_k^n}\left(\frac{\b^{\gammab}\partial^{\gammab}P\left(\hat{\y}(i)\right)}{\gammab !}\right)^{u_{k,\gammab}}
\right]\\
&=&\sum_{\u =(\u_1,\dots, \u_d) \in V(\alphab)} \frac{\alphab !}{\prod_{k=1}^d \u_k !}\prod_{k=1}^d \left[
\prod_{\gammab \in \Delta_k^n}\left(\frac{\partial^{\gammab}P\left(\hat{\y}(i)\right)}{\gammab !}\right)^{u_{k,\gammab}}\right]\\
&&\qquad\qquad\qquad\times\b^{\sum_{k=1}^d\sum_{\gammab \in \Delta_k^n} u_{k,\gammab} \gammab}\\
&=& \sum_{\u  \in V(\alphab)} P_{\alphab, \u}^i \left(\hat{\y}(i)\right) ~\b^{g(\u)}.
\end{eqnarray*}
Substituting this into \eqref{yi+bbeta-N}, and combining with \eqref{decoforsing}, 
we obtain the assertion of Proposition \ref{ypaqas}.\qed
\bigskip


The second ingredient for the proof of Theorem \ref{mainbissection3}
is the following result by E. Friedman and A. Pereira, which is a
clever use of a Raabe-type formula. 
Applying Lemma \ref{lem_mahler} to $P_\a, Q_\a$ we see that
both $Y(P_\a,Q_\a;s)$ and $Z(P_\a,Q_\a;s)$ converge absolutely in $\left\{\Re(s) >\frac{n+q}{d}\right\}$ and have meromorphic continuation to $\C$ with at most simple poles located in the set 
 $$\left\{\frac{n+q -k}{d}; ~k\in \N_0\right\}\setminus (-\N_0).$$
 
\begin{lemma}\label{zpaqas}{\rm (Friedman and Pereira \cite{fredman}, Proposition 2.2)}\\
For any $N\in \N_0$,
 $$\a \in [0,\infty)^n\mapsto Y(P_\a,Q_\a; -N) \quad {\mbox { and }} \quad \a \in [0,\infty)^n\mapsto Z(P_\a,Q_\a; -N)$$
 are  polynomials in $\a$.
If we write 
 $$Y(P_\a,Q_\a; -N)=\sum_{\alphab\in\N_0^n} c_{\alphab} \a^{\alphab}=\sum_{\alphab} c_{\alphab} \prod_{i=1}^n a_i^{\alpha_i},$$
then 
$$
Z(P,Q; -N)=\sum_{\alphab\in\N_0^n} c_{\alphab} \prod_{i=1}^n B_{\alpha_i},
$$
where the $B_k$ are the classical Bernoulli numbers defined by \eqref{def_Bernoulli}. 
\end{lemma}
We need in fact the following version of Lemma \ref{zpaqas}.
\begin{lemma}\label{zpaqasbis}
If we write 
 $$Y(P_\a,Q_\a; -N)=\sum_{\alphab\in\N_0^n} d_{\alphab} \prod_{i=1}^n (1+a_i)^{\alpha_i},$$
then 
$$
Z(P,Q; -N)=\sum_{\alphab\in\N_0^n} d_{\alphab} \prod_{i=1}^n {\widetilde{B}}_{\alpha_i}.
$$
\end{lemma}
{\bf Proof of Lemma \ref{zpaqasbis}:}\\
First we remark that 
\begin{eqnarray*}
Y(P_\a,Q_\a; -N)&=&\sum_{\alphab\in\N_0^n} d_{\alphab} \prod_{i=1}^n (1+a_i)^{\alpha_i}
=\sum_{\alphab\in\N_0^n } d_{\alphab} \sum_{\k \in \N_0^n \atop \k \leq \alphab}  \left(\prod_{i=1}^n{\alpha_i \choose k_i}\prod_{i=1}^n a_i^{k_i}\right)\\
&=&\sum_{\k\in\N_0^n } \left(\sum_{ \alphab \geq \k}  d_{\alphab} \prod_{i=1}^n{\alpha_i \choose k_i}\right) \prod_{i=1}^n a_i^{k_i}.
\end{eqnarray*}
Lemma \ref{zpaqas} then implies that 
\begin{eqnarray}\label{proof_lemma2}
Z(P,Q; -N)&=&\sum_{\k\in\N_0^n } \left(\sum_{ \alphab \geq \k}  d_{\alphab} \prod_{i=1}^n{\alpha_i \choose k_i}\right) \prod_{i=1}^n B_{k_i}
= \sum_{\alphab\in\N_0^n }  d_{\alphab} \sum_{ \k \leq \alphab } \left( \prod_{i=1}^n{\alpha_i \choose k_i} \prod_{i=1}^n B_{k_i}\right)\nonumber\\
&=& \sum_{\alphab\in\N_0^n }  d_{\alphab} \prod_{i=1}^n \left(\sum_{ k_i=0}^{\alpha_i}  {\alpha_i \choose k_i} B_{k_i}\right).
\end{eqnarray}
Since it is well known that
\begin{align}\label{bernoulliformula}
\sum_{k=0}^\alpha {\alpha \choose k} B_{k}=(-1)^{\alpha}B_{\alpha}=\widetilde{B}_{\alpha},
\end{align}
we see that the right-hand side of \eqref{proof_lemma2} is equal to
$$
\sum_{\alphab\in\N_0^n} d_{\alphab} \prod_{i=1}^n {\widetilde{B}}_{\alpha_i}.
$$
This ends the proof of Lemma  \ref{zpaqasbis}.\qed

Theorem \ref{mainbissection3}} follows  clearly from Proposition \ref{ypaqas} and Lemma \ref{zpaqasbis}.\qed

\section{Proof of Theorem \ref{mainsection3}}\label{sec-proof-firstmaintheorem}

Consider for any $j=1,\dots, n$ a polynomial $ P_j\in \R[X_1,\dots, X_j]$ in $j$ variables.
Assume that the assumptions (\ref{growthcondition}), (\ref{forsimplicity}) hold. 
Assume also that for all $j=1,\dots n-1$ the polynomial $P_j$ satisfies the assumption $(H_0S)$ and that the polynomial 
$P_n$ is elliptic and homogeneous of degree $d\geq 1$.

As mentioned in the Introduction, 
$$\zeta_n(\s; \P)= \sum_{ m_1, \dots, m_n \geq 1} \frac{1}{\prod_{j=1}^n P_j(m_1, \dots, m_j)^{s_j}}$$
converges  absolutely in the domain 
$\mathcal D_n(\P)$ (defined by (\ref{def_DnP})),
and has the meromorphic continuation to the whole space $\C^n$. 
Fix $\delta_1,\dots, \delta_n>0$ such that (\ref{growthconditionbis}) holds.  Fix $\Nb=(N_1,\dots, N_n) \in \N_0^n$ and set 
$$\sigma_0:=\frac{1}{\delta_n} \max \left\{ n+1-j+ \sum_{i=j}^n \delta_i N_i \mid j\in \{1,\dots,n\}\right\}.$$
Define for $t\in\C$
$$\psi_{\Nb, \P}(t):=\zeta_n (-\Nb+t \e_n; \P) = \zeta_n\left((-N_1,\dots, -N_{n-1}, -N_n+t); \P\right).$$
It follows from the above that 
$\psi_{\Nb, \P}(t)$ converges absolutely in  $\{\Re(t)>\sigma_0\}$ and has a meromorphic continuation to the whole complex plane $\C$.

Consider also the zeta function 
$$t\mapsto Z(P_n, \prod_{j=1}^n P_j^{N_j}; t):=\sum_{ m_1, \dots, m_n \geq 1} \frac{\prod_{j=1}^n P_j(m_1, \dots, m_j)^{N_j}}{ P_n(m_1, \dots, m_n)^{t}}$$
in one variable $t$ and set $\nu_1:=n+q_{\Nb}= n+\sum_{j=1}^n N_j \deg (P_j)$.
Lemma \ref{lem_mahler} implies that $Z(P_n, \prod_{j=1}^n P_j^{N_j}; t)$ converges absolutely in $\{\Re(t)>\nu_1/\deg P_n\}$ and has the meromorphic 
continuation to the whole $t$-plane with at most simple poles located on the set 
$$\mathcal P(\Nb; \P):=
\left\{\frac{\nu_1 -k}{\deg P_n}\mid k\in \N_0\right\}\setminus (-\N_0).$$

We have for any $t\in \C$ such that $\Re(t) >\max \{\sigma_0, \nu_1/\deg P_n\}$,
\begin{eqnarray*}
\psi_{\Nb, \P}(t)&=&
\sum_{ m_1, \dots, m_n \geq 1} \frac{1}{\left(\prod_{j=1}^{n-1} P_j(m_1, \dots, m_j)^{-N_j}\right) P_n(m_1, \dots, m_n)^{-N_n+t} }\\
&=&\sum_{ m_1, \dots, m_n \geq 1} \frac{\prod_{j=1}^n P_j(m_1, \dots, m_j)^{N_j}}{ P_n(m_1, \dots, m_n)^{t}} 
=Z(P_n, \prod_{j=1}^n P_j^{N_j}; t),
\end{eqnarray*}
and hence 
we deduce by analytic continuation that 
$$\psi_{\Nb, \P}(t) = Z(P_n, \prod_{j=1}^n P_j^{N_j}; t) \quad \mbox{for all}\;\; t\in \C\setminus \mathcal P(\Nb; \P).$$
It follows that there exists $\eta =\eta( \Nb, \P) >0$ such that 
$$\psi_{\Nb, \P}(t) = Z(P_n, \prod_{j=1}^n P_j^{N_j}; t) \quad \mbox{for all}\;\; t\in D^*(0,\eta)=\{ t\in \C;~0\leq|t|<\eta\},$$
because $t=0$ is not included in $\mathcal P(\Nb; \P)$.
Theorem  \ref{mainbissection3} implies then that 
$$\zeta_n^{\e_n}(-\Nb; \P) :=\lim_{t\rightarrow 0}\zeta_n\left((-N_1,\dots, -N_{n-1}, -N_n+t); \P\right)=\lim_{t\rightarrow 0}\psi_{\Nb, \P}(t)$$
exists and is given  by 
$$\zeta_n^{\e_n}(-\Nb; \P) =Z(P_n, \prod_{j=1}^n P_j^{N_j}; 0).$$
We conclude the assertion of Theorem \ref{mainsection3} by using the expression of $Z(P_n, \prod_{j=1}^n P_j^{N_j}; 0)$ given by Theorem \ref{mainbissection3}.  \qed

\section{Examples of transcendental values}\label{exampletranscendence}

We already mentioned in Section \ref{sec-firstmainextension} (Corollary 
\ref{irrationalityCoro}) that some special values of the power-sum zeta-function
\eqref{NEZMzetadef} are transcendental.    The same phenomenon is therefore expected
for more general zeta-function $\zeta_n(\s; \P) $.
It is surely difficult to prove the transcendency in general,
but at least we can present some examples.
The first example, essentially due to Pierrette Cassou-Nogu\`es,
gives a transcendental values at the origin.
\medskip

{\bf Example 1:}
Let $\P=(P_1,P_2,P_3,P_4)$ with
$P_1 \in \R[X_1], P_2 \in {\mathbb R}[X_1, X_2]$, $P_3 \in \R[X_1, X_2, X_3]$ 
satisfying   $(H_0S)$
and $P_4=X_1^3+X_2^3+X_3^3+X_4^3$.

Consider $\zeta_4(\s; \P)$.
We have for $\Nb=\zerob =(0,0,0,0)$:
\begin{equation}\label{4-1}
\zeta_4^{\e_4}(\zerob; \P)  :=\lim_{t\rightarrow 0}\zeta_4(t\e_4; \P)=Z(P_4,0),
\end{equation}
where 
$Z(P_4,s)$ is the meromorphic continuation of the one variable Dirichlet series 
$$s\mapsto \sum_{ m_1, \dots, m_4 \geq 1} \frac{1}{P_4(m_1, \dots, m_4)^{s}}.$$
Pierrette Cassou-Nogu\`es \cite{CN1} proved that 
$$Z(P_4,0)=
-\frac{8}{9}B_4\Gamma \left(1/3\right)^3+B_1^4=
\frac{4}{135} \Gamma \left(1/3\right)^3 +\frac{1}{16}.$$
By using the already quoted result of Chudnovsky that $\Gamma \left(1/3\right)$ is a transcendental number, it follows that 
$Z(P_4,0)$ is a transcendental number.
Therefore from (\ref{4-1}) we deduce that 
$\zeta_4^{\e_4}(\zerob; \P)$ 
is also a transcendental number.
\medskip

Next we consider the individual period integral \eqref{period_def}, and give an example
whose value is transcendental.

{\bf Example 2:}
Let $\P=(P_1,P_2,P_3)$ with $P_1=X_1$, $P_2=X_1+X_2$ and 
$P_3=X_1^2+2X_1X_2+X_2^2+X_3^2$.   Then $n=3$ and $d=\deg P_3=2$.
We consider the value of $\zeta_3(\s; \P)$ at the point $-\Nb$, where
$\Nb=(1,1,0)$.   Then $Q_{\Nb}=X_1(X_1+X_2)$ and $q_{\Nb}=\deg Q_{\Nb}=2$.

We choose one of the period integrals on the right-hand side of \eqref{explicitvalue} and
evaluate its value.   The vectors $\alphab=(\alpha_1,\alpha_2)$ and 
$\betab=(\beta_1,\beta_2,\beta_3)$ should satisfy 
$|\betab|\leq q_{\Nb}=2$ and $\alpha_1+2\alpha_2+|\betab|=q_{\Nb}+n=5$.
We choose $\betab=(2,0,0)$ and $\alphab=(1,1)$.

Next we write $\Delta_1^3=\{\gamma(11),\gamma(12),\gamma(13)\}$ and
$\Delta_2^3=\{\gamma(21),\ldots,\gamma(26)\}$, where 
\begin{align*}
&\gamma(11)=(1,0,0), \gamma(12)=(0,1,0),\gamma(13)=(0,0,1),\\
&\gamma(21)=(2,0,0), \gamma(22)=(0,2,0),\gamma(23)=(0,0,2),\\
&\gamma(24)=(1,1,0), \gamma(25)=(1,0,1),\gamma(26)=(0,1,1).
\end{align*}
Then $V(\alphab)$ is the set of vectors $\u=(\u_1,\u_2)$ with
$\u_1=(u_{1,\gamma(11)},u_{1,\gamma(12)},\textcolor{red}{u_{1,\gamma(13)}})$, 
$\u_2=(u_{2,\gamma(21)},\ldots,u_{2,\gamma(26)})$, and 
$$
\sum_{j=1}^3 u_{1,\gamma(1j)}=\alpha_1=1,\quad
\sum_{j=1}^6 u_{2,\gamma(2j)}=\alpha_2=1.
$$
We choose $\u$ where $u_{1,\gamma(11)}=u_{2,\gamma(21)}=1$ and all other 
components are 0.

Now we compute the integral $K_3(P_3;Q_{\Nb};0;\alphab,\u,\betab)$.
First, we find that
$$
P_3(\hat{\y}(3))^{-|\alphab|}=(y_1^2+2y_1y_2+y_2^2+1)^{-2}.
$$
Secondly, since $Q_{\Nb}(\hat{\y}(3))=y_1(y_1+y_2)$, we have
$$
\partial^{\betab}Q_{\Nb}(\hat{\y}(3))=\frac{\partial^2}{\partial y_1^2}
y_1(y_1+y_2)=2.
$$
Thirdly,
\begin{align*}
P_{\alphab,\u}^3(\hat{\y}(3))=\frac{1}{2!}\frac{\partial}{\partial y_1}
P_3(\hat{\y}(3))\frac{\partial^2}{\partial y_1^2}P_3(\hat{\y}(3))=2(y_1+y_2).
\end{align*}
Therefore
$$
K_3(P_3;Q_{\Nb};0;\alphab,\u,\betab)=4\int_0^1\int_0^1((y_1+y_2)^2+1)^{-2}
(y_1+y_2)dy_1 dy_2.
$$ 
Putting $y_1+y_2=z$ we see that the right-hand side is
\begin{align*}
&=4\int_0^1 dy_2\int_{y_2}^{y_2+1}\frac{z}{(z^2+1)^2}dz
=2\int_0^1\left(\frac{1}{y_2^2+1}-\frac{1}{(y_2+1)^2+1}\right)dy_2\\
&=2(2\arctan 1-\arctan 2)=2\left(\frac{\pi}{2}-\arctan 2\right)
=2\arctan\frac{1}{2}.
\end{align*}
If $\xi:=\arctan(1/2)$ is an algebraic number, then $\tan\xi$ is a transcendental
number by the classical Lindemann-Weierstrass theorem.    However
$\tan\xi=\tan(\arctan(1/2))=1/2$ is obviously not transcendental.
Therefore $\xi$ is, and hence $K_3(P_3;Q_{\Nb};0;\alphab,\u,\betab)$ is a
transcendental number.

\section{Some relations among Bernoulli numbers}\label{bernoulli}

As we mentioned in the introduction, in our previous paper \cite{GMZV}, we proved an 
explicit formula for $\zeta_n^{\thetab} (-\Nb; \gammab; \b)$
(the generalized Euler-Zagier type) in terms of $\Nb$, $\thetab$ 
and Bernoulli numbers $B_n$.
Since our $\zeta_n(\s; \P) $ includes $\zeta_n(\s; \gammab; \b)$ as
special examples, Theorem \ref{mainsection3} in the present paper also gives 
an explicit expression for
$\zeta_n^{\thetab} (-\Nb; \gammab; \b)$.
The arguments in those two papers are methodologically similar (both based on a kind of
Raabe-type formula), but not exactly the same, and consequently, the two expressions
are different.     Comparing these two expressions, we find some non-trivial relations
among Bernoulli numbers.
\medskip

{\bf Example 3:}
The case of the Riemann zeta-function $\zeta(s)$.
Applying Theorem 1 in \cite{GMZV}, for any $N\in\N_0$ we obtain
\begin{align}\label{B1}
\zeta(-N)=-\frac{1}{N+1}B_{N+1}-\sum_{\alpha=0}^N\binom{N}{\alpha}\frac{1}{N+1-\alpha}
B_{\alpha}.
\end{align}

On the other hand, let apply Theorem \ref{mainsection3} to $\zeta(s)$.
Then $n=1$, $\P=P_1$ with $P_1(X_1)=X_1$, $d=1$, $Q_N=X_1^N$, $q_N=N$.
Therefore the possible values of $\beta$ in the formula \eqref{explicitvalue} are
$\beta=0,1,\ldots,N$, and for each $\beta$, $\alpha$ is determined by
$\alpha+\beta=N+1$.    Also $V(\alpha)=\{\alpha\}$, $g(u)=\alpha$, $\hat{\y}=y_1$ and
$\hat{\y}(1)=(1)$.   Since $\partial P_1(X_1)=1$ and 
$$\partial^{\beta}Q_N(X_1)=(N)_{\beta}X_1^{N-\beta}$$ 
(where $(N)_{\beta}=1$ if $\beta=0$ and $(N)_{\beta}=N(N-1)\cdots(N-\beta+1)$ if $\beta>0$), we have
$P_1(\hat{\y}(1))=1$, $(P_1)_{\alpha,u}^1(\hat{\y}(1))=1$, and
$\partial^{\beta}Q(\hat{\y}(1))=(N)_{\beta}$, and hence the period
$K_1(P_1;Q_N;0;\alpha,u,\beta)=(N)_{\beta}$.   Therefore
\begin{align}\label{B2}
\zeta(-N)&=
\sum_{\beta=0}^N\frac{(-1)^{N+1-\beta}(N-\beta)!}{\beta!(N+1-\beta)!}
(N)_{\beta}\widetilde{B}_{N+1}\notag\\
&=\frac{B_{N+1}}{N+1}+\sum_{\beta=1}^N \frac{(-1)^{\beta}}{\beta!}
N(N-1)\cdots(N-\beta+2)B_{N+1}\notag\\
&=\frac{1}{N+1}\sum_{\beta=0}^N (-1)^{\beta}\binom{N+1}{\beta}B_{N+1}\notag\\
&=\frac{1}{N+1}\left(\sum_{\beta=0}^{N+1} (-1)^{\beta}\binom{N+1}{\beta}-(-1)^{N+1}\right)B_{N+1}=\frac{(-1)^N}{N+1}B_{N+1},
\end{align}
which coincides with the well-known classical expression of $\zeta(-N)$.   Comparing 
\eqref{B1} and \eqref{B2}, we obtain
$$
-B_{N+1}-\sum_{\alpha=0}^N\binom{N}{\alpha}\frac{N+1}{N+1-\alpha}B_{\alpha}
=(-1)^N B_{N+1},
$$
which implies the known but non-trivial formula \eqref{bernoulliformula}.
\medskip

{\bf Example 4:}
The case of the Euler double zeta-function
$$
\zeta_2(\s;\P)=\sum_{m_1=1}^{\infty}\sum_{m_2=1}^{\infty}
\frac{1}{m_1^{s_1}(m_1+m_2)^{s_2}},
$$
that is, $\s=(s_1,s_2)$, $\P=(P_1,P_2)$, where $P_1(X_1)=X_1$, $P_2(X_1,X_2)=X_1+X_2$.
First we apply Theorem 1 in \cite{GMZV} to evaluate the value at $(s_1,s_2)=-\Nb=(-N_1,-N_2)$
($N_1,N_2\in\N_0$).    In \cite{GMZV} we consider the general directional limit
$\lim_{t\to 0}(-N_1+t\theta_1,-N_2+t\theta_2)$, but in the present paper we only consider
the case $\theta_1=0$, that is $\zeta_2^{\e_2}(-\Nb;\P)$ under the notation of the
present paper.    Therefore we just state the case $\theta_1=0$. 
With notations of \cite{GMZV}, we have $\zeta_2(\s;\P)=\zeta_2(\s; \gammab; \b)$ with $\gammab=(1,1)$ and $\b=(0,1)$. 
The consequence of Theorem 1 in \cite{GMZV} is then
\begin{align}\label{B3}
&\zeta_2^{\e_2}(-\Nb;\P)\notag\\
&=\sum_{l=0}^{N_1}\sum_{\substack{k_1\geq l,k_2\geq 0\\k_1+k_2\leq N_1+N_2+2}}
\frac{(N_1+N_2+2-l)!}{(k_1-l)!k_2!(N_1+N_2+2-k_1-k_2)!}\notag\\
&\qquad\times\binom{N_1}{l}{\binom{N_1+N_2+1-l}{N_2}}^{-1}
\frac{(-1)^{N_1+3-l}B_{k_1}B_{k_2}}{(N_1+N_2+2-l)(l-N_1-1)}\notag\\
&\;\;+\sum_{l=0}^{N_1}\sum_{\substack{k_1\geq l,k_2\geq 0\\k_1+k_2\leq N_2+1+l}}
\frac{(N_2+1)!}{(k_1-l)!k_2!(N_2+1+l-k_1-k_2)!}\notag\\
&\qquad\times\binom{N_1}{l}
\frac{B_{k_1}B_{k_2}}{(N_2+1)(N_1+1-l)}\notag\\
&\;\;+\sum_{l_1=0}^{N_1}\sum_{l_2=0}^{N_2}
\sum_{\substack{k_1\geq l_1,k_2\geq 0\\k_1+k_2\leq l_1+l_2}}
\frac{l_2!}{(k_1-l_1)!k_2!(l_1+l_2-k_1-k_2)!}\notag\\
&\qquad\times\binom{N_1}{l_1}\binom{N_2}{l_2}
\frac{B_{k_1}B_{k_2}}{(N_1+N_2+2-l_1-l_2)(N_2+1-l_2)}.
\end{align}

Next we apply Theorem \ref{mainsection3} to $\zeta_2(\s;\P)$.    Then $n=2$, $d=1$,
$Q_{\Nb}=X_1^{N_1}(X_1+X_2)^{N_2}$, $q_{\Nb}=N_1+N_2$, so
$\betab=(\beta_1,\beta_2)$ should satisfy $\beta_1+\beta_2\leq N_1+N_2$ and
$\alpha=N_1+N_2+2-\beta_1-\beta_2$.
It follows that $\Delta_1^2=\{\gamma(11),\gamma(12)\}$ with $\gamma(11)=(1,0)$ and
$\gamma(12)=(0,1)$.
The set $V(\alpha)$ can be parametrized as
$$
V(\alpha)=\{\u_l=(\alpha-l,l)\;|\;0\leq l\leq \alpha\}.
$$
Then $g(\u_l)=(g_1(\u_l),g_2(\u_l))$ with $g_1(\u_l)=\alpha-l$ and
$g_2(\u_l)=l$, and we obtain
\begin{align}\label{B4}
&\zeta_2^{\e_2}(-\Nb;\P)\notag\\
&=\sum_{\substack{\beta_1,\beta_2\geq 0 \\ \beta_1+\beta_2\leq N_1+N_2}}
\sum_{l=0}^{N_1+N_2+2-\beta_1-\beta_2}
\frac{(-1)^{N_1+N_2+2-\beta_1-\beta_2}(N_1+N_2+1-\beta_1-\beta_2)!}
{(N_1+N_2+2-\beta_1-\beta_2)!\beta_1!\beta_2!}\notag\\
&\qquad\times\left(\sum_{i=1}^2 K_i(P_2;Q_{\Nb};0;\alpha,\u_l,\betab)\right)
\widetilde{B}_{N_1+N_2+2-\beta_2-l}\widetilde{B}_{\beta_2+l}.
\end{align}

We compute $K_i(P_2;Q_{\Nb};0;\alpha,\u_l,\betab)$.    First we see that
if $\beta_2>N_2$ then $\partial^{\betab}Q_{\Nb}=0$, and if $\beta_2\leq N_2$ then
\begin{align*}
&\partial^{\betab}Q_{\Nb}=(N_2)_{\beta_2}\sum_{j=\max\{0,\beta_1-N_1\}}
^{\min\{\beta_1,N_2-\beta_2\}}\binom{\beta_1}{j}(N_1)_{\beta_1-j}(N_2-\beta_2)_{j}
X_1^{N_1-\beta_1+j}(X_1+X_2)^{N_2-\beta_2-j}.
\end{align*}
Next, since 
$\partial^{\gamma(11)}P_2=\partial^{\gamma(12)}P_2=1$, we have
$$
(P_2)_{\alpha,\u_l}^i(\hat{\y}(i))=\frac{\alpha!}{(\alpha-l)!l!}=
\binom{N_1+N_2+2-\beta_1-\beta_2}{l}.
$$
Therefore
\begin{align}\label{B5}
&K_1(P_2;Q_{\Nb};0;\alpha,\u_l,\betab)\notag\\
&=\binom{N_1+N_2+2-\beta_1-\beta_2}{l}(N_2)_{\beta_2}\sum_{j=\max\{0,\beta_1-N_1\}}
^{\min\{\beta_1,N_2-\beta_2\}}\binom{\beta_1}{j}(N_1)_{\beta_1-j}(N_2-\beta_2)_{j}\notag\\
&\times\int_0^1 (1+y_2)^{-\alpha+N_2-\beta_2-j}dy_2,
\end{align}
and the last integral is
$$
=\int_0^1(1+y_2)^{-N_1-2+\beta_1-j}dy_2=\frac{1-2^{-N_1-1+\beta_1-j}}{N_1+1+j-\beta_1}.
$$
We also obtain an expression for $K_2(P_2;Q_{\Nb};0;\alpha,\u_l,\betab)$, almost the
same as \eqref{B5}, the only difference is that the corresponding integral factor is
$$
\int_0^1(1+y_1)^{-N_1-2+\beta_1-j}y_1^{N_1-\beta_1+j}dy_1
=\frac{2^{-N_1-1+\beta_1-j}}{N_1+1+j-\beta_1}.
$$
Therefore
\begin{align*}
&K_1(P_2;Q_{\Nb};0;\alpha,\u_l,\betab)+K_2(P_2;Q_{\Nb};0;\alpha,\u_l,\betab)\\
&\;=\binom{N_1+N_2+2-\beta_1-\beta_2}{l}(N_2)_{\beta_2}\sum_{j=\max\{0,\beta_1-N_1\}}
^{\min\{\beta_1,N_2-\beta_2\}}\binom{\beta_1}{j}(N_1)_{\beta_1-j}(N_2-\beta_2)_{j}\notag\\
&\qquad\qquad\times \frac{1}{N_1+1+j-\beta_1}\\
&\;=\binom{N_1+N_2+2-\beta_1-\beta_2}{l}(N_2)_{\beta_2}\sum_{j=\max\{0,\beta_1-N_1\}}
^{\min\{\beta_1,N_2-\beta_2\}}\binom{\beta_1}{j}(N_1)_{\beta_1-j-1}(N_2-\beta_2)_{j},
\end{align*}
which with \eqref{B4} implies
\begin{align}\label{B6}
&\zeta_2^{\e_2}(-\Nb;\P)\notag\\
&=\sum_{\substack{\beta_1,\beta_2\geq 0 \\ \beta_1+\beta_2\leq N_1+N_2}}
\sum_{l=0}^{N_1+N_2+2-\beta_1-\beta_2}
\frac{(-1)^{N_1+N_2+2-\beta_1-\beta_2}(N_1+N_2+1-\beta_1-\beta_2)!}
{(N_1+N_2+2-\beta_1-\beta_2)!\beta_1!\beta_2!}\notag\\
&\qquad\times\binom{N_1+N_2+2-\beta_1-\beta_2}{l}(N_2)_{\beta_2}\notag\\
&\qquad\times\sum_{j=\max\{0,\beta_1-N_1\}}
^{\min\{\beta_1,N_2-\beta_2\}}\binom{\beta_1}{j}(N_1)_{\beta_1-j-1}(N_2-\beta_2)_{j}
\widetilde{B}_{N_1+N_2+2-\beta_2-l}\widetilde{B}_{\beta_2+l}.
\end{align}

Comparing \eqref{B3} and \eqref{B6}, we obtain

\begin{prop}\label{bernoulliidentity}
The right-hand side of \eqref{B3} is equal to the right-hand side of \eqref{B6}.
This identity gives a non-trivial relation among Bernoulli numbers.
\end{prop}

Using the known data of the values of Bernoulli numbers, we can check, for instance, 
that both \eqref{B3} and \eqref{B6} gives
$\zeta_2^{\e_2}(\mathbf{0};\P)=5/12$, which agrees with the known ``reverse'' value.
However the authors do not know whether the identity obtained in Proposition 
\ref{bernoulliidentity} is essentially new, or can be deduced from known formulas.

For each multiple zeta-function $\zeta_n(\s; \P) $, we can argue as above and can obtain
certain (more and more complicated) identity among Bernoulli numbers.


\verb??\\
{\bf Driss Essouabri}\\
Univ. Lyon,\\
UJM-Saint-Etienne,\\
CNRS, Institut Camille Jordan UMR 5208,\\
Facult\'e des Sciences et Techniques,\\       
23 rue du Docteur Paul Michelon,\\
F-42023, Saint-Etienne, France.\\
{\it E-mail address}: driss.essouabri@univ-st-etienne.fr

\medskip

\verb??\\
{\bf Kohji Matsumoto}\\
Graduate School of Mathematics\\
Nagoya, University\\
Furo-cho, Chikusa-ku\\
Nagoya 464-8602, Japan\\
{\it E-mail address}: kohjimat@math.nagoya-u.ac.jp

 \end{document}